\documentclass[leqno]{article} 

\usepackage[utf8]{inputenc} 
\usepackage{amsmath} 
\usepackage{amsthm} 
\usepackage[upint]{stix2} 
\usepackage{inconsolata} 
\usepackage[margin=1.5in]{geometry}
\usepackage{mathtools}
\usepackage{enumerate} 
\usepackage{xcolor}
\usepackage[colorlinks,allcolors=blue]{hyperref} 
\usepackage[nameinlink]{cleveref} 
\usepackage[url=false,eprint=false,style=alphabetic,maxbibnames=99]{biblatex} 

\usepackage[final]{microtype}
\usepackage{todonotes} 

\swapnumbers 

\usepackage{bookmark}
\bookmarksetup{
  numbered, 
  open,
}

\usepackage{tikz-cd}

\usepackage{verbatim}

\numberwithin{equation}{subsection} 
\newtheorem{theorem}[equation]{Theorem} 
\newtheorem{lemma}[equation]{Lemma}
\newtheorem{corollary}[equation]{Corollary}

\newtheorem{proposition}[equation]{Proposition}

\theoremstyle{definition} 
\newtheorem{definition}[equation]{Definition}
\newtheorem{construction}[equation]{Construction}

\newtheorem{example}[equation]{Example}

\newtheorem{remark}[equation]{Remark}

\newtheorem{notation}[equation]{Notation}

\renewcommand{\in}{\smallin}

\newcommand{\categ}{\mathsf}

\newcommand{\E}[1]{\mathbb{E}_{#1}}

\newcommand{\Prlk}{\categ{Pr^L_k}}
\newcommand{\Prl}{\categ{Pr^L}}
\DeclareMathOperator{\LMod}{LMod}
\DeclareMathOperator{\BiMod}{BiMod}
\DeclareMathOperator{\RMod}{RMod}
\DeclareMathOperator{\LModk}{LMod_k}
\DeclareMathOperator{\QCoh}{QCoh}
\newcommand{\nPrlk}[1]{\categ{Pr^{L, {#1}}_k}}
\newcommand{\bignPrlk}[1]{\categ{Pr^{L, {#1} \kern 1pt \wedge}_k}}

\DeclareMathOperator{\cts}{cts}

\newcommand{\nPrlpoint}[2]{\categ{Pr^{L, {#1}}_{{#2}, \ast}}}

\newcommand{\nPrlkpoint}[1]{\categ{Pr^{L, {#1}}_{k, \ast}}}
\newcommand{\ctsCat}{\categ{Cat}_{\cts}}
\newcommand{\bigctsCat}{\widehat{\ctsCat}}
\newcommand{\nCatk}[1]{\categ{Cat^{#1}_k}}
\newcommand{\nPrlA}[1]{\categ{Pr^{L, {#1}}_A}}

\DeclareMathOperator{\Fun}{Fun}
\DeclareMathOperator{\id}{id}
\DeclareMathOperator{\op}{op}
\DeclareMathOperator{\Mod}{Mod}

\DeclareMathOperator{\Spc}{\categ{Spc}}

\DeclareMathOperator{\bigSpc}{\widehat{\Spc}}
\DeclareMathOperator{\Hom}{\mathcal{H}\kern -2pt om}
\DeclareMathOperator{\Spec}{Spec}
\DeclareMathOperator{\Maps}{Maps}
\DeclareMathOperator{\ObjDef}{ObjDef}
\DeclareMathOperator{\fmpObjDef}{\widehat{\ObjDef}}
\DeclareMathOperator{\SimDef}{SimDef}
\DeclareMathOperator{\sDef}{sDef}
\DeclareMathOperator{\fmpSimDef}{\widehat{\SimDef}}
\DeclareMathOperator{\CatDef}{CatDef}
\DeclareMathOperator{\fmpCatDef}{\widehat{\CatDef}}
\newcommand{\ObjDefMap}{\beta^{\operatorname{obj}}}
\newcommand{\SimDefMap}{\beta^{\operatorname{sim}}}
\DeclareMathOperator{\AlgDef}{AlgDef}
\DeclareMathOperator{\fmpAlgDef}{\widehat{\AlgDef}}

\newcommand{\nMonDef}[1]{\operatorname{\mathbb{E}_{#1}{\text -}MonDef}}
\newcommand{\nMonCat}[1]{\operatorname{\mathbb{E}_{#1}{\text -}MonCat}}

\DeclareMathOperator{\Alg}{Alg}

\DeclareMathOperator{\aug}{aug}
\DeclareMathOperator{\End}{End}
\newcommand{\naugAlg}[1]{\operatorname{{Alg^{({#1}), aug}_k}}}
\newcommand{\nAlgslice}[1]{\operatorname{{Alg^{({#1})}_{k/k}}}}
\newcommand{\naugAlgslice}[2]{\operatorname{{Alg^{({#1}), aug}_{k/{#2}}}}}

\newcommand{\nsmAlg}[1]{\operatorname{{Alg^{({#1}), sm}_k}}}
\newcommand{\nsmAlgslice}[2]{\operatorname{{Alg^{({#1}), sm}_{k/{#2}}}}}
\newcommand{\nAlg}[1]{\operatorname{{Alg^{({#1})}_k}}}

\newcommand{\C}{\categ{C}}
\newcommand{\D}{\categ{D}}
\newcommand{\M}{\categ{M}}

\DeclareMathOperator{\LM}{LM}
\DeclareMathOperator{\RM}{RM}
\DeclareMathOperator{\Assoc}{Assoc}

\DeclareMathOperator{\alg}{alg}
\DeclareMathOperator{\LModalg}{LMod^{\alg}}
\DeclareMathOperator{\LModaug}{LMod^{\aug}}

\DeclareMathOperator{\RModalg}{RMod^{\alg}}
\DeclareMathOperator{\RModaug}{RMod^{\aug}}

\DeclareMathOperator{\Arr}{Arr}
\DeclareMathOperator{\Triv}{Triv}

\DeclareMathOperator{\LCatn}{LCat^n}
\DeclareMathOperator{\LCatnstar}{LCat^{n}_{\ast}}
\DeclareMathOperator{\LCatnstarcocart}{LCat^{n, \cocart}_{\ast}}

\DeclareMathOperator{\RCatn}{RCat^n}
\DeclareMathOperator{\RCatnstar}{RCat^{n}_{\ast}}
\DeclareMathOperator{\RCatnstartriv}{RCat^{n,triv}_{\ast}}

\newcommand{\Dual}{\mathcal{D}}
\newcommand{\lindual}[1]{{#1}^{\vee}}

\DeclareMathOperator{\Z}{\xi}
\DeclareMathOperator{\Free}{Free}
\DeclareMathOperator{\AlgLMC}{\Alg_{/ \LM}(\C^\otimes)}
\DeclareMathOperator{\AlgRMC}{\Alg_{/ \RM}(\C^\otimes)}
\DeclareMathOperator{\AlgAssocC}{\Alg_{\Assoc / \LM}(\C^\otimes)}

\DeclareMathOperator{\AlgAssocRC}
{\Alg_{\Assoc / \RM}(\C^\otimes)}

\DeclareMathOperator{\ev}{ev}
\DeclareMathOperator{\Deform}{Deform}
\DeclareMathOperator{\cocart}{cocart}
\DeclareMathOperator{\pair}{\mathcal{M}}
\DeclareMathOperator{\fib}{Fib}
\DeclareMathOperator{\barconst}{Bar}
\DeclareMathOperator{\LinFun}{LinFun}
\title{Deformations of objects in $n$-categories}
\author{Dennis Chen}
\date{April 1, 2023}

\begin{document}

\maketitle
\begin{abstract}
    In this paper, we prove that the deformation theory of an object in an $n$-category is controlled by the its $n$-fold endomorphism algebra. This recovers Lurie's results on deforming objects and categories. We also generalize a previous result by Blanc et al. (\cite{Blanc}) on deforming a category and an object simultaneously to the case of $n$-categories. 
\end{abstract}
\section{Introduction}
In algebraic geometry, there is a notion of deforming various objects over local Artinian algebras. For example, take a scheme $X$ and a quasicoherent module $M$ over it. A deformation of $M$ over the dual numbers $k[\epsilon]$ is the data of a quasicoherent module $M_{\epsilon}$ over 
\[
X_{\epsilon} := \Spec{k[\epsilon]} \times X
\] 
whose pullback along the inclusion $X \to X_{\epsilon}$ gives $M$. If one requires $M, M_\epsilon$ to be locally free, then $M_\epsilon$ is characterized wholly by its gluing data, which in this case is captured in the first cohomology group of $\End(M)$. Hence locally free deformations over $k[\epsilon]$ of a locally free module is characterized by classes in $H^1(\End(M))$.

Lurie generalizes these examples to the case of deforming an object in a category (\cite[Section~5.2]{DAGX}, \cite[Section~16.5]{SAG}). Using the framework of formal moduli problems, he shows that the $\E{1}$-formal moduli problem associated to deforming an object can be characterized by its algebra of endomorphisms:
\[
\fmpObjDef_M \simeq \Maps_{\naugAlg{1}}(\Dual({-}), k \oplus \End_\C(M)).
\] 

There's also the classical notion of deforming a category and relating it to its Hochschild complex, as explained in \cite{Kontsevich, Seidel, KL}, which is important for example in the study of Mirror symmetry and Fukaya categories.

Lurie reformulates this result in the context of infinity categories (\cite[Section~5.3]{DAGX}, \cite[Section~16.6]{SAG}):
\[
\fmpCatDef_{\C} \simeq \Maps_{\naugAlg{2}}(\Dual({-}), k \oplus \xi(\C))
\] where here $\xi(\C)$ denotes the derived center of $\C$, which can be calculated via the Hochschild complex of $\C$.

In this paper we follow Lurie's arguments to generalize and unify his results of deforming an object in a category and deforming a category in $\Prl$. Namely, given any $k$-linear $n$-category $\C$ and an object $M \in \C$, we construct a functor $\ObjDef_M$ and show that the $\E{n}$ algebra characterizing the formal moduli problem (\cite[Definition~12.1.3.1]{SAG}) associated to $\ObjDef_M$ is the $n$-iterated endormorphism space of $M$, or the center of $M$. We recover Lurie's results for object deformations by taking $\C$ to be a $1$-category, and we recover his results for category deformations by taking $\C$ to be $\Prlk$. 

More precisely, let $k$ be a field, $\Prlk$ be the $(\infty, 1)$-category of presentable $k$-linear categories and $k$-linear colimit preserving functors. It has a monoidal structure given by the $k$-linear tensor of categories. In other words, $\Prlk := \LMod_{\LModk}(\Prl)$, the category of presentable categories with $\Mod_k$ action. Then inductively, we define $\nPrlk{n}$ as the $(\infty, 1)$-category of presentable linear categories tensored over $\nPrlk{n-1}$ (see \ref{def:nprlk}), the objects of which we call "$k$-linear $n$-categories". In a similar vein, one can define $\nPrlA{n}$ for an $\E{n+1}$-algebra $A$. 

We define a version of object deformations for a $k$-linear $n$-category. That is, given an $n$-category $\C \in \nPrlk{n}$ and an object $M \in \C$, we define a functor
\[
\ObjDef_M: \nsmAlg{n} \to \bigSpc
\]
from small $\E{n}$-algebras to large spaces.
The functor is intuitively given by the formula
\[
\ObjDef_M(A) := \LMod_A(\C)\times_\C \{M\}.
\]
This functor will have an associated formal moduli problem which is characterized by an augmented $\E{n}$-algebra: the $n$-fold endomorphism algebra of $M$ (thought of as a nonunital algebra). More precisely, we have the following definition:

\begin{definition}[$n$-fold Endomorphism object] \label{def:end}
Given $M \in \C$.
Let $\End^1_\C(M) = \Hom_\C(M, M)$. This has a clear basepoint $\id_M$
Inductively we can define 
\[
\End^{n+1}_\C(M) := \Hom_{\End^n_\C(M)}(\id^n_M, \id^n_M),
\]
where $\id^n_M$ (or sometimes $1^n_M$) is the identity of $\End^n_\C(M)$, with a new basepoint given by the identity $\id^{n+1}_M \in \Hom_{\End^n_\C(M)}(\id, \id)$. When the context is clear, we may drop the $n$ and $M$ from $\id^n_M$.

We may also use $\End^0_\C(M)$ to denote $\C$ where the basepoint is $M$, which is an alternate base case for this induction. Here $\Hom$ denotes the internal hom, see (\ref{notation:general}).
\end{definition}

Our first main result is:
\begin{theorem}\label{theorem:objdef1}
The formal moduli problem associated to $\ObjDef_M$ is equivalent to 
\[
\Maps_{\naugAlg{n}}(\Dual^n({-}), k \oplus \End_\C^n(M)),
\]
where $\Dual^n$ is the $\E{n}$-Koszul duality functor (\cite[Section~5.2.5]{HA}).
\end{theorem}

This directly generalizes previous results: using $n = 1$ we get exactly the classical result for deforming objects in categories \cite[Section~5.2]{DAGX}, \cite[Section~16.5]{SAG}. Using $n=2$ and $\C = \nPrlk{2}$ and letting $M$ be a given category in $\C$, we get Lurie's result for deforming categories \cite[Section~5.3]{DAGX}, \cite[Section 16.6]{SAG}. The proof is given in \ref{theorem:objdef2}.

We next consider the problem of deforming an object and $n$-category simultaneously. We follow Blanc, Katzarkov, and Pandit (\cite[Section~4]{Blanc}) and Lurie (\cite[Remark~16.0.0.3]{SAG}), who previously considered the case of $n = 1$. 

More precisely, let
$\nPrlpoint{n}{A}$ for the category of presentable pointed $A$-linear $n$-categories:
\[
\nPrlpoint{n}{A} := \categ{Pr^{L, {n}}_{{A}, \LMod^{n}_{A} /}}
\]

Given an $n$-category $\C$ and an object $M \in \C$, we can define a simultaneous deformation functor as follows: given a small $\E{n+1}$-algebra $A$, we let
\[
\SimDef_{(\C, M)}(A) := \nPrlpoint{n}{A} \times_{\nPrlkpoint{n}} \{(\C, M)\}
\]
where the map 
\[
\nPrlpoint{n}{A} \to \nPrlkpoint{n}
\]
is using the augmentation $A \to k$, and for consistency with the our other section, we let 
\[
\nPrlpoint{n}{A} := \LMod_{\LMod^n_A} (\nPrlkpoint{n})
\]
be the category of {\it left} $A$-modules in $\nPrlkpoint{n}$ (as opposed to \cite{Blanc} which uses right modules). We show that the formal completion of this functor is characterized by the nonunital $\E{n+1}$-algebra 
\begin{equation}\label{eq:Z}
\Z(\C, E) := \fib(\Z(\C) \to \Z(M))
\end{equation}
where $\Z(\C)$ is the center of $\C$, $\Z(M)$ is the center of $M$, and the fiber is taken at $0 \in \Z(M)$. Explicitly, we can let $\Z(\C) := \End^{n+1}_{\nPrlk{n}}(\C)$,  $\Z(M):=\End^n_{\C}(M)$ and the map between them is given by evaluation at $M$. 

Our second main result is:
\begin{theorem}\label{theorem:simdef1}
There is an equivalence of formal moduli problems:
\[
\fmpSimDef_{(\C, M)} \to \Maps_{\naugAlg{n+1}}(\Dual^{n+1}({-}), k \oplus \Z(\C, M)).
\]
\end{theorem}
The proof is given in (\ref{theorem:simdef2}). For example, using $n = 1$, the center of $\C$ is represented by $\End_{\End(C)}(1_\C)$---in other words natural transformations from $1_C$ to itself---and the center of $M$ is represented by $\End_{\C}(M)$. The map $\Z(\C) \to \Z(M)$ is given by evaluation of the natural transformation at $M$. This recovers Proposition 4.7 of \cite{Blanc}. 

These deformation problems are related to the deformation problem of an $\E{n}$-monoidal category:
Given an $\E{n}$-monoidal category $\D$, its deformations can be identified with deformations of the pointed category $(\LMod^n_{\D}, \LMod^{n-1}_{\D})$. These ideas are discussed in section \ref{section:3.5}. This uses the fully faithful embedding of $\E{n}$-monoidal categories into $n$-pointed categories (categories with an object together) via the rule
\[
\D^\otimes \mapsto (\LMod^n_{\D}, \LMod^{n-1}_{\D}).
\]
Hence one can study the deformation theory of $\D$ by studying the deformations of the pair $(\LMod^n_{\D}, \LMod^{n-1}_{\D})$. This also recovers the deformations of $\E{n}$ algebras (which can be thought of as single-object $n$-categories $k$-cells are trivial for $k<n$). To\"en, in theorem 5.1 and 5.2 of \cite{Toen}, 
also relates deformations of $\E{n}$-monoidal categories to $\E{n+1}$-Hochschild cochains as defined in \cite{Francis}. 

The deformation theory of $\E{n}$-monoidal categories is incredibly important for various theories of quantization. In section 2 of \cite{Toen}, To\"en explains the connection between different variations of quantization---namely quantum groups, skein algebras, and Donaldson-Thomas invariants---to deformations of (monoidal) categories. For example, To\"en relates quantum groups (see \cite{Drinfeld}) to deforming the category of sheaves of the moduli space $\operatorname{Bun}_G(\ast)$ of $G$-bundles on the point. 

\subsection{Acknowledgements}
I'd like to thank my advisor, David Nadler, from whose guidance I have benefited enormously. In addition, I'd like to thank Germ\'an Stefanich for many insightful conversations and ideas. Lastly, this work was partially supported by NSF RTG grant DMS-1646385

\subsection{Set theoretic issues} 
For this section, let's hypothesize for now an increasing sequence of universes $U_0, U_1, \dots$ We let "small" mean $U_0$-small and "large" mean $U_1$-small. We will only need two universes $U_0$ and $U_1$. 

To solve the set theoretic issues of even defining $\nPrlk{n}$ (and $\nPrlA{n}$ by analogy), we follow \cite{German}. There are two solutions. 

First we can define $\nCatk{1}$ as the category of $k$-linear categories with cocontinuous $k$-linear functors between them. Then $\nCatk{2}$ to be the $U_2$-small category of all $\Prlk$-linear $U_1$-small categories with $\Prlk$-linear cocontinuous functors. We can continue the induction, producing larger and larger categories $\nCatk{n}$ which is $U_n$-small. Notice there is no presentability here.

The other idea is to only use two universes, one small $U_0$ and one large $U_1$. Then one can define $\nPrlk{n}$ inductively, following chapter 12 of \cite{German}. 

\begin{definition}[Presentable $k$-linear $n$-categories]\label{def:nprlk}
Let's define $\nPrlk{n}$ and $\bignPrlk{n}$ inductively:
For $n=0$, we define
\[
\nPrlk{0} = \bignPrlk{0} := \Mod_k.
\]
Next we inductively define:
\[
\bignPrlk{n} := \Mod_{\nPrlk{n-1}}(\bigctsCat).
\]
In other words, $\bignPrlk{n}$ is the category of $\nPrlk{n-1}$-modules in the large category of cocomplete categories and cocontinuous functors.
Finally we can define $\nPrlk{n}$ to be the full subcategory of $U_0$-compact objects of $\bignPrlk{n}$. We will also denote this category by $\LMod^{n+1}_k$, see (\ref{remark:nprlk_and_lmod}).
\end{definition}

\begin{remark}
Notice that $\nPrlk{1}$ agrees with the usual definition of $\Prlk$ as the category of presentable $k$-linear categories with cocontinuous $k$-linear functors between them.
\end{remark}

\begin{remark}\label{remark:nprlk}
We will make use of $\bignPrlk{n}$ mainly because like in \cite[Remark~8.4.3]{German}, we don't know if the hom objects for $\C \in \nPrlk{n}$ are presentable. They may only exist in $\bignPrlk{n-1}$. However it is true that for any category $\D \in \bignPrlk{n}$, its hom objects are in $\bignPrlk{n-1}$.
\end{remark}

Despite which approach we take, our $\ObjDef$ fmp is perhaps large in general, in contrast to \cite[Section~16.5, 16.6]{SAG}. The point is that our $n$-categories as defined could have large $n$-fold endomorphism objects, unlike the category deformation and object deformation problems that Lurie considered. The author doesn't know whether these presentable $n$-categories have presentable hom objects or not. 

However, if our given $n$-category $\C$ had small $n$-fold endomorphism objects, then we have the following easily using our main result \ref{theorem:objdef1}:

\begin{proposition}
Given an $n$-category $\C \in \nPrlk{n}$ (or $\nCatk{n}$) with small $n$-fold endomorphism objects. Then the formal moduli problem associated to $\ObjDef_{M \in \C}$ is a functor that lands in $\Spc$, the category of $U_0$-small spaces. 
\end{proposition}

In this paper, we will by default use the second method of restricting to presentable categories for concreteness, but the arguments don't really differ regardless of which method one chooses. 

\subsection{Conventions}
Our notational conventions are listed here. First, unless otherwise mentioned, we are working over a field $k$ and all mentions of $k$-linear objects are infinity categorical. For example, "finite dimensional vector space" will mean a compact object in the infinity category $\Mod_k$ of $k$-modules. 

We will also occasionally use the abbrevation "fmp" for "formal moduli problem". 

By default, categories of algebras will be large due to $\ObjDef$ being large. Notice that $\nsmAlg{n}$ is still a small category due to the finiteness conditions placed on small algebras.
\subsubsection{Notation} \label{notation:general}

Here is some basic notation and conventions. We will have more notation later which will be introduced as needed.
\begin{itemize}
    \item $\C, \D$ describe categories in $\nPrlk{n}$ (\ref{def:nprlk}). We also just call these $n$-categories for short. 
    \item Similarly pairs $(\C, E), (\D, F)$ denote objects in $\nPrlkpoint{n}$.
    \item $\Spc$ denotes the category of small spaces. $\bigSpc$ denotes the category of $U_1$-small spaces, or "large" spaces. 
    \item $\Maps_{\C}(x, y) \in \bigSpc$ is the space of maps between $x$ and $y$.
    \item $\Hom_{\C}(x, y)$ denotes the hom object in $\nPrlk{n-1}$, as given by the right adjoint to the tensor action on $\C \in \nPrlk{n}$. 
    \item $\Z(M)$ denotes the center of $M \in \C$ which can be calculated via $\End^n_{\C}(M)$.
    \item $\Z(\C, E)$ denotes the center of $(\C, E) \in \nPrlkpoint{n}$ which can be calculated via 
    \[
    \Z(\C, E) := \fib(\Z(\C) \to \Z(M)).
    \]
    \item $\Alg^n$ denotes the large category of $\E{n}$-algebras (or more precisely, $n$-fold iterated algebras), following \cite{HA}. $\nAlg{n}$ denotes the large category of $\E{n}$-algebras over $k$.
    \begin{itemize}
        \item $\naugAlg{n}$ denotes the (large) category of augmented $\E{n}$ algebras over $k$.
        \item $\nsmAlg{n}$ denotes the category of artinian/small augmented $\E{n}$-algebras. Notice that this is always a small category due to the definition of small algebras. 
    \end{itemize}
    \item $\LMod(\C)$ denotes $\AlgLMC$ for a category $\C \in \nPrlk{n}$, following  Definition 4.2.1.13 of \cite{HA}.
    \begin{itemize}
        \item The cocartesian fibration 
        \begin{equation}\label{equation:lmodfib}
            \LMod(\C) \to \Alg(\nPrlk{n-1})
        \end{equation}
        denotes the cocartesian fibration
        \[
            \AlgLMC \to \AlgAssocC,
        \] again following \cite{HA} and using that $\AlgAssocC \simeq \Alg(\nPrlk{n-1})$ for $\C \in \nPrlk{n}$.
    \end{itemize}
    \item $\LMod^n_R$ denotes the $n$-fold left modules of an $\E{n}$-algebra $R$ (see \ref{definition:lmod}). We let $\LMod^0_R$ denote $R$, and $\LMod^{-1}_R$ denote the unit $1_R \in R$. For $R = k$, notice that $\LMod^n_k$ agrees with $\nPrlk{n-1}$, see (\ref{remark:nprlk_and_lmod}).
    \item $M_R$ denotes $\LMod_R^{n-1} \otimes M$ for $M \in \C \in \nPrlk{n}$ and $R \in \nAlg{n-1}$.
    \item $\C_A$ denotes $\LMod^n_A \otimes \C \in \nPrlA{n}$.
    \item $(\C_A, E_A)$ denotes $(\LMod^n_A \otimes \C, \LMod^{n-1}_A \otimes E) \in \nPrlpoint{n}{A}$.
    \item $\widehat{F}$: given a functor $F: \nsmAlg{n} \to \bigSpc$, we let $\widehat{F}$ denote its fmp completion (see \cite[Remark~1.1.17]{DAGX}, \cite[Remark~12.1.3.5]{SAG}). In other words, it is the formal moduli problem associated to $F$. For example, $\fmpObjDef_M$ is the fmp completion of $\ObjDef_M$.
\end{itemize}

\section{Deformations of objects}
Throughout this section, we assume we are given $\C \in \nPrlk{n}$ and an object $M \in \C$. We would like to show that the deformations of $M$ in the $n$-category $\C$ is characterized by the $\E{n}$-algebra $\End_\C^n(M)$.
The argument follows four steps, generally following Lurie's ideas in \cite[Section~5.2, 5.3]{DAGX} or \cite[Section~16.5, 16.6]{SAG}:
\begin{enumerate}
    \item Construct the functor $\ObjDef_M$ characterizing deformations of $M \in \C$.
    \item Prove $\ObjDef_M$ is $n$-proximate.
    \item Construct the comparison map \[\ObjDefMap: \ObjDef_M \to \Maps_{\naugAlg{n}}(\Dual^n({-}), k \oplus \End_\C^n(M))\]
    \item Prove that $\ObjDefMap$ is an equivalence. 
\end{enumerate}
\subsection{Constructing the functor \texorpdfstring{$\ObjDef_M$}{ObjDefM}}

First let's define how to take $n$-fold modules for an $\E{n}$ algebra. We have functors
\[
\Alg(\nPrlk{n}) \to \nPrlk{n+1}
\]
via $\categ{D} \mapsto \LMod_\categ{D}(\nPrlk{n})$. We also get induced functors by applying $\nAlg{p}$ to get
\[
\Alg^{p+1}(\nPrlk{n}) \to \Alg^p(\nPrlk{n+1}).
\]
Composing these functors, we can define:
\begin{definition}[Iterated left modules]\label{definition:lmod}
Let $\LMod^n$ denote the composite functor given by 
\begin{equation} 
\LMod^n: \nAlg{n} \to \nAlg{n-1}(\Prlk) \to \nAlg{n-2}(\nPrlk{2}) \dots \to \Alg(\nPrlk{n-1}) \to \nPrlk{n}.
\end{equation}
We denote the evaluation of this functor on $R$ by $\LMod^n_R$. 
We will also use the variant functor 
\begin{equation} 
\LMod^{n-1}: \nAlg{n} \to \nAlg{n-1}(\Prlk) \to \nAlg{n-2}(\nPrlk{2}) \dots \to \Alg(\nPrlk{n-1}).
\end{equation}
which ends one step early as compared to $\LMod^n$ above. 
\end{definition}

\begin{remark}
Clearly there are variants where one can take iterated right modules, or even switch between taking left and right modules. Notice that if $R$ is an $\E{\infty}$-algebra, then these constructions are all equivalent, and we may denote the category by $\Mod^n_R$.
\end{remark}
\begin{remark} \label{remark:nprlk_and_lmod}
For $R = k$, notice that we get $\LMod^n_k \simeq \nPrlk{n-1}$ for $n \geq 1$. We can show this by induction: 
For $n = 1$, clearly $\LMod^1_k = \Mod_k =: \nPrlk{0}$. Next, if the result holds true for $n=m-1$, then we notice that
\begin{align*}
    \LMod^m_k & := \LMod^1_{\LMod^{m-1}_k} \\
    & \simeq \LMod_{\nPrlk{m-2}}(\nPrlk{m-1}) \\
    & \simeq \nPrlk{m-1}
\end{align*}
since every object of $\nPrlk{m-1}$ is already an $\nPrlk{m-2}$-module, by definition (\ref{def:nprlk}). Since $k$ is an $\E{\infty}$-algebra, we may also denote $\LMod^n_k$ by $\Mod^n_k$. 
\end{remark}
Recall that we define $\LMod(\C):= \AlgLMC$ (\ref{equation:lmodfib}), hence we have a cocartesian fibration
\[
\LMod(\C) \to \Alg(\nPrlk{n-1}).
\]
We pull back along $\LMod^{n-1}$ to get the left modules whose action is given by an $\E{n}$-algebra.
\begin{definition}[Left Modules]\label{def:lmodstar}
Let $\LModalg(\C)$ to be the following pullback:
\[
\begin{tikzcd}
\LModalg(\C) \ar[r] \ar[d] \arrow[dr, phantom, "\scalebox{1.5}{$\lrcorner$}" , very near start, color=black]
& \LMod(\C) \ar[d]
\\ \nAlg{n} \ar[r] & \Alg(\nPrlk{n-1})
\end{tikzcd}
\] Here the left vertical map is the cocartesian fibration $\LMod(\C) \to \AlgAssocC$ and the lower horizontal map is the $n$-fold $\LMod$ functor.

We can also further pull back along $\naugAlg{n} \to \nAlg{k} \to \Alg(\nPrlk{n-1})$ to get 
\[
\begin{tikzcd}
\LModaug(C) \ar[r] \ar[d] \arrow[dr, phantom, "\scalebox{1.5}{$\lrcorner$}" , very near start, color=black]
& \LMod(\C) \ar[d]
\\ \naugAlg{n} \ar[r] & \Alg(\nPrlk{n-1})
\end{tikzcd}
\]
\end{definition}

Intuitively, objects of $\LModalg(C), \LModaug(\C)$ consists of triples $(A, E, \eta)$ where $E \in \C$, $A \in \nAlg{n}$ (or $\naugAlg{n}$), and $\eta$ is a left action of $\LMod^{n-1}_A$ on $E$. 

\begin{remark}\label{remark:rmod}
Dually, we can use right modules instead by replacing the right vertical leg by
\[
\AlgRMC \to \AlgAssocRC
\] 
where we are using the fact that since $\nPrlk{n}$ is symmetric monoidal, so we can choose to either use left or right modules to define $\C^{\otimes}$.

In other words, we can pull back along the cocartesian fibration
\[
\RMod(\C) \to \Alg(\nPrlk{n-1}).
\]
So we can analoguously define $\RModalg(\C)$ and $\RModaug(\C)$. This has objects $(A, E, \eta)$ where $E \in \C$, $A \in \nAlg{n}$ (or $\naugAlg{n}$), and $\eta$ is a {\it right} action of $\LMod^{n-1}_A$ on $E$.
\end{remark}

Now we are finally ready to construct our functor. To do this, we construct first the associated fibration, then use straightening/unstraightening to get the functor we need.
Recall we have a cocartesian fibration $\LModalg(C) \to \naugAlg{n}$. Let $\LModalg(C)^{\cocart}$ be the subcategory whose morphisms are the cocartesian arrows of this map. Then
\begin{equation}\label{equation:lmodstar}
\LModalg(C)^{\cocart} \to \naugAlg{n}    
\end{equation}
is a left fibration. 
Our given object $M$ has a natural $\LMod^{n-1}_k$ action. This gives us an object $(k, M) \in \LModalg(C)$. 
\begin{definition}[Deformation fibration]\label{definition:objfib}
Let $\Deform[\C, M]$ be the slice of $\LModalg(C)^{\cocart}$ over the object $(k, M)$. In other words:
\[
\Deform[\C, M] := (\LModalg(\C)^{\cocart})_{/(k, M)}.
\]
We have an induced left fibration $\Deform[\C, M] \to \naugAlg{n}$ by taking slices of \ref{equation:lmodstar}:
\[
\Deform[\C, M] := \LModalg(C)^{\cocart}_{/(k, M)} \longrightarrow \nAlgslice{n} \simeq \naugAlg{n}.
\]
\end{definition}
\begin{construction}[$\ObjDef$]\label{construction:objdef}
$\Deform[\C, M] \to \naugAlg{n}$ classifies a functor $\widetilde{\Deform}[\C, M]: \naugAlg{n} \to \bigSpc$. Here $\bigSpc$ is the category of not-necessarily $U_0$-small spaces. Finally, by restricting to small algebras, we get the functor we wanted: 
\[
\ObjDef_M: \nsmAlg{n} \to \bigSpc.
\]
\end{construction}

\begin{remark}
Notice that given an algebra $A \in \naugAlg{n}$, $\ObjDef_M(A) \simeq \LMod^n_A(\C) \times_\C \{M\}$, where the map $\LMod^n_A(\C) \to \LMod^n_k(\C) \simeq \C$ is given by the augmentation map $A \to k$. Here $\LMod^n_k \simeq \nPrlk{n-1}$ (\ref{remark:nprlk_and_lmod}), so since $\C$ is $\nPrlk{n-1}$-linear , every object is a $\nPrlk{n-1}$-module. 
\end{remark}
\subsection{Proving \texorpdfstring{$\ObjDef_M$}{ObjDefM} is \texorpdfstring{$n$}{n}-proximate} \label{section:2.2}

Let's begin with a generalization of fully faithfulness. Here we let $\C(x, y) := \Hom_{\C}(x, y)$ for brevity.
For the following definiton and proposition, we'll need to use non-presentable $N$-categories (in other words, we consider categories in $\bignPrlk{n}$) because we'll induct by taking repeated hom spaces (\ref{remark:nprlk}).

\begin{definition}[$n$-fully faithful] \label{def:ff}
Given a functor $F; \C \to \D$ of $N$-categories. Then we say
\begin{itemize}
    \item $F$ is $0$-fully faithful if $F$ is an equivalence.
    \item $F$ is $n$-fully faithful if for all $x, y \in \C$. the induced functor $\C(x, y) \to \D(F(x), F(y))$ is $(n-1)$-fully faithful.
\end{itemize}
This is an inductive definition for $0 \leq n \leq N$. Notice that the $n=1$ case agrees with our usual notion of fully faithfulness.
\end{definition}

Now we prove that $\ObjDef_M$ is an $n$-proximate fmp. We recall the following result:
\begin{proposition}[\cite{SAG} Prop 16.2.1.1]\label{prop:SAG}
Let $k$ be an $\E{2}$-ring and $\D$ a $k$-linear category. Suppose we are given a pullback:
\[
\begin{tikzcd}
    A \ar[r] \ar[d] \arrow[dr, phantom, "\scalebox{1.5}{$\lrcorner$}" , very near start, color=black] 
    & A' \ar[d]
\\  B \ar[r] & B'
\end{tikzcd}
\] in $\Alg_k$. Then the induced functor $\LMod_A(\D) \to \LMod_{A'}(\D) \times_{\LMod_{B'}(\D)}  \LMod_B(\D)$ is fully faithful. 
\end{proposition}
\begin{remark}\label{remark:ff}
    Using the hypotheses of the above proposition (\ref{prop:SAG}), let $M \in \LMod_A(D)$, and let $M_R := R \otimes_A M$ for any $A$-algebra $R$. We also denote $M_A = M$ to match the above notation. 
    Then, as explained in the proof of \cite[Prop~16.2.1.1]{SAG}, the conclusion of the above proposition (\ref{prop:SAG}) is equivalent to the unit map 
    \[
    M_A \longrightarrow M_B \times_{M_{B'}} M_{A'}
    \]
    being an equivalence for any $M \in \LMod_A(D)$. This is an easy application of the result that left adjoints are fully faithful if and only if the unit map is an equivalence. 
\end{remark}
Using this result, we'll prove the following result by induction:
\begin{proposition}\label{prop:ff}
Given a pullback in $\nAlg{n}$,
\[
\begin{tikzcd}
    A \ar[r] \ar[d] \arrow[dr, phantom, "\scalebox{1.5}{$\lrcorner$}" , very near start, color=black] 
    & A' \ar[d]
\\  B \ar[r] & B'
\end{tikzcd}
\]
and $\C \in \bignPrlk{n}$, with $M \in \LMod^{n-1}_A \otimes C$. Given $R \in \nAlg{n}$ under $A$, write $M_R$ for 
\[
\LMod^{n-1}_R \otimes M.
\]
We also use $M_A$ for $M$ to match the above notation.
Then the comparison map 
\[
M_A \longrightarrow M_B \times_{M_{B'}} M_{A'}
\] is representably $(n-1)$-fully faithful.
\end{proposition}

Here "representably" $n$-fully faithful means after taking homs out from any object $\C(X, {-})$, the result is a $n$-fully faithful functor.

\begin{proof}
The case $n=1$ is explained in remark (\ref{remark:ff}). This is the base for our induction.

For the inductive step, we prove it for $n > 1$, assuming it is done for $n - 1$. Then given our pullback square along with $M \in \C$, we are trying to show that 
\[
M_A \longrightarrow M_B \times_{M_{B'}} M_{A'}
\] 
is representably $(n-1)$-fully faithful. This means that given any $X \in \C$, we must show 
\[
\C(X, M_A) \longrightarrow \C(X, M_B) \times_{C(X, M_{B'})} C(X, M_{A'})
\]
is $(n-1)$-fully faithful. Let $X_R$ denote $\LMod^{n-1}_R \otimes X$ as with $M$, and let $\C_R$ denote $\LMod^n_R \otimes \C$. Then by the extension of scalars adjunction, we see that 
\[
\C(X, M_R) \simeq \C_R(X_R, M_R),
\]
so our above map is equivalent to 
\[
\C_A(X_A, M_A) \longrightarrow \C_B(X_B, M_B) \times_{C_{B'}(X_{B'}, M_{B'})} C_{A'}(X_{A'}, M_{A'}).
\]
To show this map is $(n-1)$-fully faithful, we take any two objects $P_A, Q_A \in \C_A(X_A, M_A)$ and we try to show the induced map 
\[
\C_A(X_A, M_A)(P_A, Q_A) \longrightarrow \C_B(X_B, M_B)(P_B, Q_B) \times_{C_{B'}(X_{B'}, M_{B'})(P_{B'}, Q_{B'})} C_{A'}(X_{A'}, M_{A'})(P_{A'}, Q_{A'})
\] is $(n-2)$-fully faithful (using the inductive definition of $(n-1)$-fully faithfulness). Again we are using $Q_R := \LMod^{n-2}_R \otimes Q_A$, and the same with $P_R$ for $P_A$.

This last map, once again by the extension of scalars adjunction, can be identified with 
\[
\C_A(X_A, M_A)(P_A, Q_A) \longrightarrow \C_A(X_A, M_A)(P_A, Q_B) \times_{C_{A}(X_{A}, M_{A})(P_{A}, Q_{B'})} C_{A}(X_{A}, M_{A})(P_{A}, Q_{A'}).
\]

But since by induction, we assume $Q_A$ is representably $(n-2)$-fully faithful in $\C(X_A, M_A)$. This implies that the comparison morphism---the image under $\C_A(X_A, M_A)(P_A, {-})$ of $Q_A \to Q_B \times_{Q_{B'}} Q_{A'}$---is indeed $(n-2)$-fully faithful, as desired.
\end{proof}

Using this we prove 
\begin{proposition}
$\ObjDef_M$ is an $n$-proximate fmp.
\end{proposition}
\begin{proof}
We seek to prove that given a pullback:
\[
\begin{tikzcd}
    A \ar[r] \ar[d] \arrow[dr, phantom, "\scalebox{1.5}{$\lrcorner$}" , very near start, color=black] 
    & A' \ar[d]
\\  B \ar[r] & B'
\end{tikzcd}
\]
in $\naugAlg{n}$, the comparison $\ObjDef_M(A) \to \ObjDef_M(A') \times_{\ObjDef_M(B')} \ObjDef_M(B)$ is $(n-2)$-truncated.

Using our above result, we know that the comparison map
\[
M_A \longrightarrow M_B \times_{M_{B'}} M_{A'}
\]
is representably $(n-1)$-fully faithful for any $M_A \in \LMod^{n-1}_A \otimes \C$.

Thus, by extension of scalars, we have that given any $X_A \in \C_A$ (borrowing notation from the last proof), we have 
\[
\C_A(X_A, M_A) \longrightarrow \C_B(X_B, M_B) \times_{C_{B'}(X_{B'}, M_{B'})} C_{A'}(X_{A'}, M_{A'}).
\]
is $(n-1)$-fully faithful. Plugging in $X_A = M_A$, we get that 
\[
\End_{\C_A}(M_A) \longrightarrow \End_{\C_B}(M_B) \times_{\End_{\C_{B'}}(M_{B'})} \End_{\C_{A'}}(M_{A'})
\]
is $(n-1)$-fully faithful. However, notice that using the basepoint $M_A$ for $\ObjDef_M(A)$, we see that $\Omega \ObjDef_M(A)$ at that basepoint can be identified with the fiber of $\End_{\C_A}(M_A)^{\simeq} \to \End_{\C}(M)^{\simeq}$. Hence, the above $(n-1)$-fully faithful map descends through fibers and taking cores, and we get that 
\[
\Omega \ObjDef_M(A) \longrightarrow \Omega \ObjDef_M(B) \times_{\Omega \ObjDef_M(B')} \Omega \ObjDef_M(A')
\]
is $(n-1)$-faithful. It is easy to show that this is equivalent to the map being $(n-3)$-truncated as we actually have a map of spaces. For example, when $n=2$, we know $1$-fully faithful maps between spaces are equivalent to $-1$-truncated inclusions. 

Hence, as we can take the base point $M_A \in \ObjDef_M(A)$ in the above argument, we can remove the loop space functor and see that the map
\[
\ObjDef_M(A) \longrightarrow \ObjDef_M(B) \times_{ \ObjDef_M(B')}  \ObjDef_M(A')
\]
must be $(n-2)$-truncated, as desired.
\end{proof}

\subsection{Constructing the comparison map \texorpdfstring{$\ObjDefMap$}{bObjDef}}
We construct the map 
\[
\ObjDefMap: \ObjDef_M \to \Maps_{\naugAlg{n}}(\Dual^n({-}), k \oplus \Z(M)).
\]

We begin by constructing a duality functor 
\[
\Dual^n: \Deform[\C, M]^{\op} \to \RModaug(\C) \times_{\C} \{M\}
\]
where $\RModaug(\C)$ is as defined in (\ref{remark:rmod}).
\begin{construction}[Duality functor $\Dual^n$]\label{constr:dual}
Let $\lambda^n: \pair^n \to \naugAlg{n} \times \naugAlg{n}$ be the pairing of categories inducing $\E{n}$-Koszul duality (\cite[Construction~5.2.5.32]{HA}). Objects of $\pair^n$ intuitively consist of two algebras $A, B \in \naugAlg{n}$ along with an augmentation of their tensor product: $A \otimes_k B \to k$. 

Let $A \otimes_k B \to k$ be an $\E{n}$-pairing between $A, B \in \nsmAlg{n}$ (so it is an object of $\pair^n$). Suppose we are given the data $(A, M_A, \eta)$ where $M_A \in \LMod_A(\C)$ and $\eta: k \otimes_A M_A \to M$ is an equivalence, so this data can be thought of as an object in $\ObjDef_M(A)$. 
Notice that 
\[
M_A \otimes \LMod^n_B \in {}_{LMod^n_A \otimes LMod^n_B} \BiMod_{\LMod^n_B}(\C)
\]
Thus we have 
\[
M \simeq \LMod^n_k \otimes_{A \otimes B} M_A \otimes \LMod^n_B \in \RMod_{\LMod^n_B}(\C)
\]
where the equivalence uses the given pairing and $\eta$. This construction gives a right $\LMod^n_B$ action on $M$.
This construction produces a functor:
\[
\Deform[\C, M] \times \pair^n \to \Deform[\C, M] \times (\RModaug(\C) \times_{\C} \{M\})
\]

This is a left representable pairing of categories, which induces a duality functor (by \cite[Construction~3.1.3]{DAGX})
\[
\Dual^n: \Deform[\C, M]^{\op} \to \RModaug(\C) \times_{\C} \{M\},
\]
as required.
\end{construction}

Now, we can easily get our compairson $\ObjDefMap$.
\begin{construction}[Comparison map $\ObjDefMap$]
Notice that $\Dual^n$ (\ref{constr:dual}) constructed above has codomain equivalent to \[
\naugAlgslice{n}{\Z(M)},
\]
where $\Z(M)$ is a center of $M$. 

We have a square:

\[
\begin{tikzcd}
\Deform_M^{\op} \ar[r, "\Dual^n"] \ar[d] & \naugAlgslice{n}{\Z(M)} \ar[d]
\\ \naugAlg{n}^{\op} \ar[r, "\Dual^n"] & \naugAlg{n}
\end{tikzcd}
\]
Here the top horizontal functor is the duality functor defined just above, the bottom functor is $\E{n}$-Koszul duality functor. The left and right vertical maps are canonical Cartesian fibrations. 

We restrict to small algebras: 
\[
\begin{tikzcd}
\Deform_M^{\op}\big|_{\nsmAlg{n}} \ar[r, "\Dual^n"] \ar[d] & \nsmAlgslice{n}{\Z(M)} \ar[d]
\\ \nsmAlg{n}^{\op} \ar[r, "\Dual^n"] & \nsmAlg{n}
\end{tikzcd}
\]
Note that the bottom morphism is an equivalence. This morphism of the vertical left fibrations gives us a comparison morphism $\ObjDefMap$ of the two induced functors
\begin{equation}\label{map:objbeta}
\ObjDefMap: \ObjDef_M \to \Maps_{\naugAlg{n}}(\Dual^n({-}), k \oplus \Z(M)),   
\end{equation}
as desired.
\end{construction}

\subsection{Proving \texorpdfstring{$\ObjDefMap$}{bObj} induces an equivalence}

We would like to show:
\begin{theorem}\label{theorem:objdef2}
Given any $n$-category $\C$ with an object $M \in \C$, the map $\ObjDefMap$ (\ref{map:objbeta}) induces an equivalence 
\[
\fmpObjDef_M \to \Maps_{\naugAlg{n}}(\Dual^n({-}), k \oplus \Z(M)).
\]
\end{theorem}

To do this, we first use Lurie's Proposition 1.2.10 in \cite{DAGX} to reduce to the cases where the input algebra is $k \oplus k[m]$ for $m > 0$, as values on these algebras determine the tangent complex in our current deformation context. 

Then we have a square:
\[
\begin{tikzcd}
\ObjDef_M(k \oplus k[m]) \ar[r, "\ObjDefMap"] \ar[d] & \Maps_{\naugAlg{n}}(\Dual^n(k \oplus k[m]), k \oplus \Z(M)) \ar[d, "\sim"]
\\ \Omega^n \ObjDef_M(k \oplus k[m + n]) \ar[r, "\Omega^n \ObjDefMap"] & \Omega^n \Maps_{\naugAlg{n}}(\Dual^n(k \oplus k[m + n]), k \oplus \Z(M))
\end{tikzcd}
\]
Since $\ObjDef_M$ is an $n$-proximate fmp, we can reduce our task to showing that the bottom map is an equivalence. 

So we've reduced our problem to proving the following:
\begin{proposition}\label{prop:main}
Let $(\C, M)$ be an $n$-category with an object.
Then the bottom leg of the square
\[
\begin{tikzcd}
\ObjDef_M(k \oplus k[m]) \ar[r, "\ObjDefMap"] \ar[d] & \Maps_{\naugAlg{n}}(\Dual^n(k \oplus k[m]), k \oplus \Z(M)) \ar[d, "\sim"]
\\ \Omega^n \ObjDef_M(k \oplus k[m + n]) \ar[r, "\Omega^n \ObjDefMap"] & \Omega^n \Maps_{\naugAlg{n}}(\Dual^n(k \oplus k[m + n]), k \oplus \Z(M))
\end{tikzcd}
\]
is an equivalence for all $m > 0$.
\end{proposition}
We start with some preliminary lemmas. First we need a lemma about functors out of $\LMod^n_R$, a Morita style result. This is just an $n$-categorical version of \cite[Theorem~4.8.4.1]{HA}, and indeed it follows from that result. 

\begin{theorem}\label{theorem:morita}
Let $R$ be an $\E{n}$-algebra where $n \geq 1$ and $\M \in \nPrlk{n}$ be an $n$-category. 
Then the composition
\[
\begin{aligned}
\Hom_{\nPrlk{{n}}}(\LMod^n_R, \M) & \subseteq \LinFun_{\Mod^n_k}(\LMod^n_R, \M) \\
& \to \Fun(\RMod_{\LMod^{n-1}_R}(\LMod^n_R), \RMod_{\LMod^{n-1}_R}(\M)) \\
& \to \RMod_{\LMod^{n-1}_R}(\M)
\end{aligned}
\]
is an equivalence. 
The second map uses the functoriality of $\RMod_{\LMod^{n-1}_R}$, and the third map is evaluation at the bimodule $R$. 
\end{theorem}

\begin{remark}
Of course by reversing left and right, there is an analoguous dual version of \ref{theorem:morita}. Note the difference between $\Hom_{\nPrlk{{n}}}$ and $\LinFun_{\Mod^n_k}$ is that while they are both $\Mod^n_k$-linear, functors in the former also have to preserve colimits. 
\end{remark}

\begin{proof}
We just use \cite[Theorem~4.8.4.1]{HA} directly to prove this one. Let $K$ contain all small simplices. If $R$ is an $\E{n}$-algebra, then $\LMod^{n-1}_R$ is an $\E{1}$-algebra in $\Mod^n_k$. $\M$ is also right tensored over $\Mod^n_k$ (as left and right modules over $\Mod^n_k$ are equivalent). So we directly apply the dual of \cite[Theorem~4.8.4.1]{HA} using the category $\C = \Mod^n_k$, the right module $\M$, and the algebra $\LMod^{n-1}_R \in \C$. This gives us exactly what we needed.
\end{proof}

The considerations in \cite[Section~4.8]{HA}, show that the categorical dual of $\LMod_S$ is $\RMod_S$. If we apply this here where $R = \LMod^{n-1}_R$, we see that the reason that 
\[
\RMod_{\LMod^{n-1}_R}
\]
shows up is because it is the categorical dual of $\LMod^n_R$. This motivates the following notation:

\begin{notation}[Duality]\label{notation:duality}
Let 
\[
\lindual{\LMod^{n}_R}
\]
denote the categorical dual of $\LMod^n_R$ in $\nPrlk{n}$. So in the case that $n \geq 1$, we get 
\[
\lindual{\LMod^{n}_R} \simeq \RMod_{\LMod^{n-1}_R}.
\]
In the case that $n = 0$, we would get the $k$-linear dual of $R$, if it exists.
\end{notation}

Now we can do some simple calculations of endormophism spaces. First one about endormorphisms of the unit object in $\LMod^n_R$. For the next few results, recall that we have the convention that $\LMod^0_R := R$ and $\LMod^{-1}_R := 1_R \in R$.
\begin{corollary} \label{cor:morita1}
Let $R$ be an $\E{n}$-algebra, where $n \geq 1$. Then: 
\begin{enumerate}
    \item The evaluation map
\[
\ev_{\LMod^{n-2}_R}: \End_{\LMod^n_R}(\LMod^{n-1}_R) \to \LMod^{n-1}_R.
\]
is an equivalence sending the identity to $\LMod^{n-2}_R$. 
    \item The composite map of evaluations
\[
\ev^m: \End^m_{\LMod^n_R}(\LMod^{n-1}_R) \to \LMod^{n-m}_R,
\]
is an equivalence that sends the identity to $\LMod^{n-m-1}_R$.
\end{enumerate}
\end{corollary}
\begin{remark}
    Notice for the case that $n=1$, this corollary gives the simple result that 
    \[
    \ev_{1_R}: \End_{\LMod_R}(R) \to R
    \] is an equivalence. 
\end{remark}
\begin{proof}
The domain of $\ev_{\LMod^{n-2}_R}$ is
\[
\Hom_{\LMod^n_R}(\LMod^{n-1}_R, \LMod^{n-1}_R)
\] by definition. 
We use the free-forgetful adjunction between $\LMod^n_R$ and $\Mod^n_k$, which gives us 
\[
\Hom_{\LMod^n_R}(\LMod^{n-1}_R, \LMod^{n-1}_R) \simeq \Hom_{\Mod^n_k}(\Mod^{n-1}_k, \LMod^{n-1}_R).
\]
Then the Morita result \ref{theorem:morita} then let's us simplify the second mapping space via evaluation to 
\[
\begin{aligned}
\ev_{\Mod^{n-2}_k}: \Hom_{\Mod^n_k}(\Mod^{n-1}_k, \LMod^{n-1}_R) & \simeq \RMod_{\Mod^{n-2}_k} (\LMod^{n-1}_R) \\
& \simeq \LMod^{n-1}_R.
\end{aligned}
\]
But this can clearly be identified with $\ev_{\LMod^{n-2}_R}$ when precomposing with the free-forgetful adjunction, thus we are done.

For the second statement, it follows from a simple induction and reduction of various endomorphism spaces using the first result. 
\end{proof}
We can secondly calculate a result about endormophisms of the augmentation module in $\LMod^n_R$.
\begin{corollary}\label{cor:morita2}
Let $R$ be an augmented $\E{n}$-algebra where $n \geq 1$. 
Then:
\begin{enumerate}
    \item The evaluation map
    \[
\ev_{\Mod^{n-2}_k}: \End_{\LMod^n_R}(\Mod^{n-1}_k) \to \lindual{\LMod^{n-1}_{\barconst{R}}}.
\]
is an equivalence which sends the identity map to the augmentation module $\Mod^{n-2}_k$. 
\item The composite map of evaluations
\[
\ev^m: \End^m_{\LMod^n_R}(\Mod^{n-1}_k) \to \lindual{\LMod^{n-m}_{\barconst^m{R}}},
\] 
is an equivalence which sends the $m$-fold identity to the augmentation module $\Mod^{n-m-1}_k$ when $m \leq n-1$.
\end{enumerate}
\end{corollary}

\begin{remark}\label{remark:morita2}
As a special case, when $m = n$ we get
\[
\ev^n: \End^n_{\LMod^n_R}(\Mod^{n-1}_k) \to {\barconst^n{R}}^{\vee} \simeq \Dual^n R
\] is an equivalence.
\end{remark}
\begin{proof}
Let's use the extension of scalars along the augmentation map $R \to k$ to identify
\[
\begin{aligned}
\Hom_{\LMod^n_R}(\Mod^{n-1}_k, \Mod^{n-1}_k) & \simeq \Hom_{\Mod^n_k}(\Mod^{n-1}_k \otimes_R \Mod^{n-1}_k, \Mod^{n-1}_k)\\
& \simeq \Hom_{\Mod^n_k}(\LMod^{n-1}_{k \otimes_R k}, \Mod^{n-1}_k)
\end{aligned}
\]
Now if $n=1$, we see the last mapping space simplifies directly to $\barconst{R}^{\vee}$. Otherwise if $n > 1$, we use the Morita result \ref{theorem:morita} to get 
\[
\begin{aligned}
\ev_{\LMod^{n-2}_{k \otimes_R k}}: \Hom_{\Mod^n_k}(\LMod^{n-1}_{k \otimes_R k}, \Mod^{n-1}_k) & \simeq \RMod_{\LMod^{n-2}_{k \otimes_R k}}(\Mod^{n-1}_k) \\
& \simeq \RMod_{\LMod^{n-2}_{\barconst{R}}},
\end{aligned}
\]
as required. Tracing the identifications, we see it indeed corresponds with $\ev_{\LMod^{n-2}_k}$, which by definition sends the identity map to the augmentation module.

For the induction we just iteratively use the first result. Notice that even if we replace $\LMod^n_R$ with it is dual $\lindual{\LMod^n_R} = \RMod_{\LMod^{n-1}_R}$, the identification above
\[
\begin{aligned}
\Hom_{\lindual{\LMod^n_R}}(\Mod^{n-1}_k, \Mod^{n-1}_k) & \simeq \Hom_{\Mod^n_k}(\Mod^{n-1}_k \otimes_R \Mod^{n-1}_k, \Mod^{n-1}_k)\\
& \simeq \Hom_{\Mod^n_k}(\LMod^{n-1}_{k \otimes_R k}, \Mod^{n-1}_k)
\end{aligned}
\]
can basically go unchanged, which is why the induction works past the second step (which requires calculating endomorphisms of the augmentation module in the category $\lindual{\LMod^{n-m}_R}=\RMod_{\LMod^{n-m-1}_R}$ when $2 \leq m \leq n-1$). 

Lastly when $m = n$, the final dual that we take is not a categorical dual but just a $k$-linear dual, thus giving us $\lindual{\barconst^n(R)}$ at the end.
\end{proof}

\begin{remark}
Note that given $R \in \naugAlg{n}$, ie an $n$-fold augmented algebra, we can take it is opposite in $n$-different ways given its $n$-commuting multiplications. If we choose the very first multiplication to take $R^{\op}$, then we get $\LMod^n_{R^{\op}} \simeq \RMod_{\LMod^{n-1}_R}$. However, regardless of which factor we take $R^{\op}$ on, we get (not canonically) equivalent algebras because any of these $\op$-funtors correspond to choosing an element on the determinant $-1$ connected component of $O(n)$, which naturally acts on $\E{n}$. Such an identification relies on a path between these two elements of $O(n)$.

We could also take the Bar construction on various multiplication levels of $R$. Notice that if we choose to take $\op$ and $\barconst$ on the first level (as we do in the above argument), we get that $(\barconst{R})^{\op} \simeq \barconst(R^{2-\op})$, where we need to take the opposite of the second multiplication of $R$. This is because in $(\barconst{R})^{\op}$, after taking Bar on the first multiplication, it is removed (or turned into a comultliplication), hence taking op afterwards affects the second multiplication of our original algebra $R$.

By using the standard calculation 
\[
\barconst(A^{\op}) \simeq \barconst(A)
\]
on $1$-algebras, we can see that
\[
\begin{aligned}
(\barconst{R})^{\op} & \simeq \barconst(R^{2-\op}) \\
& \simeq \barconst(R^{\op}) \\
& \simeq \barconst(R),
\end{aligned}
\]
which finally let's us identify $\RMod_{\LMod^{n-2}_{\barconst{R}}}$ with $\LMod^{n-1}_{\barconst{R}}$. So as long as $m \leq n-1$ in the argument above, we didn't really need the categorical duals in the endomorphism space formula in the above lemma. Using this result would have made the induction after the second step a little more symmetric-looking, however this identification relies on a choice of a path in $O(n)$ between the two different opposites that we take and isn't canonical.
\end{remark}

Next we need some results on endomorphism spaces and tensor products. These results could have been proven directly without the above corollaries, but we separated out the arguments for clarity.
\begin{lemma}[Endomorphisms and tensors 1] \label{lemma:end1}
Let $R$ be a small $\E{n}$-algebra, and $(\C, M)$ be an $n$-category with an object. Let $\C_R$ denote $\LMod^n_R(\C)$ and $M_R$ denote $LMod^{n-1}_R \otimes M$. 
Then the canonical tensoring map
\[
    i_m: \End^m_{\LMod^n_R}(\LMod^{n-1}_R) \otimes \End^m_{\C}(M) \to \End^m_{\C_R}(M_R)
\]
is an equivalence for $0 \leq m \leq n$. In this equivalence $1^{m-1}_R  \otimes 1^{m-1}_M$ goes to $1^{m-1}_{M_R}$.
\end{lemma}

\begin{proof}
We prove it by induction. For $m = 0$, the result is obvious: the map defaults to the comparison map 
\[
i_0: \LMod^n_R \otimes \C \to \C_R
\]
which is an equivalence that also sends $\LMod^{n-1}_R \otimes M$ to $M_R$, as required. 

Next let's assume it is true for $m-1$, where $1 \leq m \leq n$ and we'll prove that it is true for $m$. We can simplify our codomain through a series of steps. First, by definition we have 
\[
\End^m_{\C_R}(M_R) := \Hom_{\End^{m-1}_{\C_R}(M_R)}(1^{m-1}_{M_R}, 1^{m-1}_{M_R}).
\]

Our inductive hypothesis says $i_{m-1}$ is an equivalence. Using the functoriality of $i_{m-1}$, combined with \ref{cor:morita1}, gives us an equivalence between $\Hom_{\End^{m-1}_{\C_R}(M_R)}(1^{m-1}_{M_R}, 1^{m-1}_{M_R})$ and 
\[
\Hom_{\LMod^{n-m+1}_R \otimes \End^{m-1}_{\C}(M)}(\LMod^{n-m}_R \otimes 1^{m-1}_{M}, \LMod^{n-m}_R \otimes 1^{m-1}_{M}).
\] 
Next, using the free-forgetful adjunction, we can simplify the above mapping space to 
\[
\Hom_{\End^{m-1}_{\C}(M)}(1^{m-1}_{M}, \LMod^{n-m}_R \otimes 1^{m-1}_{M}).
\]
Now we know that $\LMod^{n-m}_R$ is dualizable for $1 \leq m \leq n$ using results on $\LMod$ in \cite[Remark~4.8.4.8]{HA}. For $n = m$, the $\LMod^{n-m}_R$ simplifies to just $R$. Since $R$ is small, it is dualizable as a $k$-module. 

Thus for all $0 \leq m \leq n$, we know that $\LMod^{n-m}_R$ is dualizable, so we can pull it out of the mapping space:
\[
\Hom_{\End^{m-1}_{\C}(M)}(1^{m-1}_{M}, \LMod^{n-m}_R \otimes 1^{m-1}_{M}) \simeq \LMod^{n-m}_R  \otimes \Hom_{\End^{m-1}_{\C}(M)}(1^{m-1}_{M}, 1^{m-1}_{M})
\]
which again using \ref{cor:morita1} we can identify with 
\[
\End^m_{\LMod^n_R}(\LMod^{n-1}_R) \otimes \End^m_{\C}(M).
\]
We can trace this comparison map backwards and we'll see it clearly sends a pair of maps to their external tensor. In other words it induces the map $i_m$. Thus we are done, we've show $i_m$ is equivalent to a composition of equivalences. 
\end{proof}

The second endomorphism result is about the augmentation module instead of the ring itself.
\begin{lemma}[Endomorphisms and tensors 2] \label{lemma:end2}
Let $R$ be a free $\E{n}$-algebra on a finite-dimensional vector space, and $(\C, M)$ be an $n$-category with an object. Let $\C_R$ denote $\LMod^n_R(\C)$ and $M_{\aug}$ denote $LMod^{n-1}_k \otimes M$. Here $\Mod^{n-1}_k$ denotes the augmentation module.  
Then the canonical tensoring map
\[
    j_m: \End^m_{\LMod^n_R}(\Mod^{n-1}_k) \otimes \End^m_{\C}(M) \to \End^m_{\C_R}(M_{\aug})
\]
is an equivalence for $0 \leq m \leq n$. In this equivalence $1^{m-1}_k \otimes 1^{m-1}_M$ goes to $1^{m-1}_{M_{\aug}}$.
\end{lemma}
\begin{proof}
We follow the last proof and start by induction. For $m=0$, the result is by definition:
\[
j_0: \LMod^n_R \otimes \C \to \C_R
\] is clearly an equivalence sending $\Mod^{n-1}_k \otimes M$ to $M_{\aug}$.

Next let's assume it is true for $m-1$. Let's simplify our codomain through a series of steps. We have
\[
\End^m_{\C_R}(M_{\aug}) := \Hom_{\End^{m-1}_{\C_R}(M_{\aug})}(1^{m-1}_{M_{\aug}}, 1^{m-1}_{M_{\aug}}).
\]
Our inductive hypothesis says $j_{m-1}$ is an equivalence, thus it is also an equivalence on mapping spaces. This result combined with our calculation  \ref{cor:morita2} gives
\[
\begin{aligned}
& \Hom_{\End^{m-1}_{\C_R}(M_{\aug})}(1^{m-1}_{M_{\aug}}, 1^{m-1}_{M_{\aug}}) \\
\simeq & \Hom_{\lindual{\LMod^{n-m+1}_{\barconst^{m-1}{R}}} \otimes \End^{m-1}_{\C}(M)}(\LMod^{n-m}_k \otimes 1^{m-1}_{M}, \LMod^{n-m}_k \otimes 1^{m-1}_{M}).
\end{aligned}
\]

Now we use the extension of scalars along $\barconst^{m-1}{R} \to k$ to simplify the right mapping space to get
\[
\Hom_{\End^{m-1}_{\C}(M)}(\LMod^{n-m}_k \otimes_{\barconst^{m-1}{R}} \LMod^{n-m}_k \otimes 1^{m-1}_{M}, 1^{m-1}_{M})
\]
which simplifies to 
\[
\Hom_{\End^{m-1}_{\C}(M)}(\LMod^{n-m}_{\barconst^m{R}} \otimes 1^{m-1}_{M}, 1^{m-1}_{M}).
\] 
Since $\LMod^{n-m}_{\barconst^m{R}}$ is dualizable (for $n = m$ we are using $R$ is free on a finite-dimensional vector space, thus $\barconst^m(R)$ is dualizable \cite[Proposition~5.2.3.15]{HA}), we can pull it out of the mapping space:
\[
\lindual{\LMod^{n-m}_{\barconst^m{R}}} \otimes \Hom_{\End^{m-1}_{\C}(M)}(1^{m-1}_{M}, 1^{m-1}_{M})
\]
But using our calculation \ref{cor:morita2}, we can clearly identify this with our domain
\[
\End^m_{\LMod^n_R}(\Mod^{n-1}_k) \otimes \End^m_{\C}(M).
\]
If you trace this calculation, you can see that this chain of equivalences is equivalent to the tensoring map $j_m$.
\end{proof}

Now we are ready to prove the proposition.
\begin{proof}[Proof of \ref{prop:main}]
we would like to show the bottom leg of \[
\begin{tikzcd}
\ObjDef_M(k \oplus k[m]) \ar[r, "\ObjDefMap"] \ar[d] & \Maps_{\naugAlg{n}}(\Dual^n(k \oplus k[m]), k \oplus \Z(M)) \ar[d, "\sim"]
\\ \Omega^n \ObjDef_M(k \oplus k[m + n]) \ar[r, "\Omega^n \ObjDefMap"] & \Omega^n \Maps_{\naugAlg{n}}(\Dual^n(k \oplus k[m + n]), k \oplus \Z(M))
\end{tikzcd}
\] is an equivalence. 

First let $R = k \oplus k[m+n]$ (which is a small algebra), so our bottom leg is now
\begin{equation} \label{map:omega_beta}
\Omega^n \ObjDefMap_R: \Omega^n \ObjDef_M(R) \to \Omega^n \Maps_{\naugAlg{n}}(\Dual^n(R), k \oplus \Z(M)).
\end{equation} 

First we identify $\ObjDefMap_R$ with the following map
\[
\begin{tikzcd}
\LMod^n_R(\C) \times_{\C} \{M\} \ar[rr, "{}_{\Dual^n R} k_R \otimes_R {-}"] && \LMod^n_{\Dual^n R}(\C) \times_{\C} \{M\}
\end{tikzcd}
\]
by using two observations: 
\begin{itemize}
    \item For the domain, $\ObjDef_M(R) \simeq \LMod^n_R(\C) \times_{\C} \{M\}$. In the pullback, the map $\LMod^n_R(\C) \to \C$ is given by the augmentation map.
    \item For the codomain, $\Maps_{\naugAlg{n}}(\Dual^n(R), k \oplus \Z(M)) \simeq \LMod^n_{\Dual^n R}(\C) \times_{\C} \{M\}$. In the pullback, the map $\LMod^n_{\Dual^n R}(\C) \to \C$ is given by the forgetful functor.
\end{itemize}
as well as just unpacking the definition of $\ObjDefMap$ (\ref{map:objbeta}). Notice here we are suppressing the $\LMod^{n-1}$ in the tensor product ${}_{\Dual^n R} k_R \otimes_R {-}$, using the fact that 
\[
\LMod^{n-1}: \Alg^{n-1} \to \nPrlk{n-1}
\]
is monoidal and fully faithful, and thinking of the $\E{n}$ Koszul duality pairing $R \otimes \Dual^n R \to k$ as giving the $\E{n-1}$ algebra $k$ two commuting central actions by $R$ and $\Dual^n R$ \cite[Proposition~5.2.5.33, Lemma~5.2.5.36]{HA}. After taking $\LMod^{n-1}$ we get exactly the correct tensoring that defines $\ObjDefMap$.

We have a triangle 
\[
\begin{tikzcd}
\LMod^n_R(\C) \ar[rd, "U"] \ar[rr, "{}_{\Dual^n R} k_R \otimes_R {-}"] && \LMod^n_{\Dual^n R}(\C) \ar[ld, "aug"] \\
& \C
\end{tikzcd}
\]
which gives $\ObjDefMap_R$ after taking fibers at $M \in \C$---for the domain along the forgetful functor $U$ to $\C$, for the codomain along the extension of scalars of the augmentation map $R \to k$ to $\C$. This triangle commutes because if you take left adjoints everywhere, you get the classical calculation that the Koszul dual of a square-zero algebra is a free algebra \cite[Proposition~4.5.6]{DAGX}

Now the horizontal map in the triangle can be identified with 
\[
\begin{tikzcd}
\LMod^n_R \otimes \C \ar[rrr, "({}_{\Dual^n R} k_R \otimes_R {-}) \otimes 1_{\C}"] &&& \LMod^n_{\Dual^n R} \otimes \C.
\end{tikzcd}
\]
For notational convenience, given any $\E{n}$ augmented algebra $B$, let $\C_B := \LMod^n_B(\C)$ and $M_B := \LMod^{n-1}_B \otimes M$, and $M_{\aug} := \LMod^{n-1}_{k, \aug} \otimes M$. For the last definition, $\LMod^{n-1}_{k, \aug}$ is the augmentation module, induced by the augmentation map $B \to k$.

We want to analyze $\Omega^n \ObjDefMap_R$ (\ref{map:omega_beta}). Our domain can be identified with the fiber of
\[
\End^n_{\C_R}(M_R) \to \End^n_{\C}(M)
\]
and the codomain can be identified with 
\[
\End^n_{\C_{\Dual^n R}}(M_{\aug}) \to \End^n_{\C}(M).
\] The asymmetry here is because in the domain, the basepoint is $M_R$ while in the codomain, the basepoint is $M_{\aug}$. 

Our map $\Omega^n \ObjDefMap$ is induced by the $n$-functoriality of $\ObjDefMap$. Indeed on $n$-cells, $\ObjDefMap$ induces a functoriality map
\[
\theta: \End^n_{\C_R}(M_R) \to \End^n_{\C_{\Dual^n R}}(M_{\aug}).
\]

As this map is induced by the tensor product $({}_{\Dual^n R} k_R \otimes_R {-}) \otimes 1_{\C}$, we have a natural square
\[
\begin{tikzcd}
\End^n_{\LMod^n_R} (\LMod_R^{n-1}) \otimes \End^n_{\C}(M) \ar[d, "L"] \ar[rr, "\alpha"] && \End^n_{\LMod^n_{\Dual^n R}} (\LMod_{k, \aug}^{n-1}) \otimes \End^n_{\C}(M) \ar[d, "R"] \\
\End^n_{\C_R}(M_R) \ar[rr, "\theta"] &&\End^n_{\C_{\Dual^n R}}(M_{\aug})
\end{tikzcd}.
\]
The left and right legs come from the from the identification $M_R := \LMod^{n-1}_R \otimes M$ and $M_{\aug} := \Mod^{n-1}_k \otimes M$. Notice $\alpha$ is induced by the functoriality of the map
\begin{equation} \label{map:alpha}
\begin{tikzcd}
\LMod^n_R \ar[rrr, "{}_{\Dual^n R} k_R \otimes_R {-}"] &&& \LMod^n_{\Dual^n R}
\end{tikzcd}
\end{equation}
on the left factor and identity on the right factor (since it is functoriality of the identity on the right factor).

We would like to show that $L, R, \alpha$ are all equivalences, which would then show $\theta$ is an equivalence. This would let us conclude $\Omega \ObjDefMap$ is an equivalence by taking fibers. 

First, showing $L$ and $R$ are equivalences is just our calculation of endomorphism spaces in \ref{lemma:end1}, and \ref{lemma:end2}.

Next we try to show $\alpha$ is an equivalence. Since the right hand factor of $\alpha$ is identity, we only have to focus on the left hand factor of $\alpha$:
\[
\alpha': \End^n_{\LMod^n_R} (\LMod_R^{n-1}) \to \End^n_{\LMod^n_{\Dual^n R}} (\LMod_{k, \aug}^{n-1}),
\] which is induced by the functoriality of
\[
\begin{tikzcd}
\LMod^n_R \ar[rrr, "{}_{\Dual^n R} k_R \otimes_R {-}"] &&& \LMod^n_{\Dual^n R}
\end{tikzcd}
\] on $n$-cells.
By Morita equivalence \ref{theorem:morita}, any $\LMod^n_k$-linear colimit preserving functor 
\[
F: \LMod^n_R \to \LMod^n_{\Dual^n R}
\]
corresponds uniquely to a bimodule structure ${}_{\Dual^n R} Q_R$ on $Q \in \LMod^n_k$ 
(we are again suppressing notation: in reality we have ${}_{\LMod^{n-1}_{\Dual^n R}} Q_{\LMod^{n-1}_R}$).

The forward direction of this equivalence takes the functor $F$ and evaluates it on $\LMod^{n-1}_R$, which gives ${}_{\Dual^n R} Q$ with its $\Dual^n R$ action. The $R$ action comes from looking at $n$-fold endomorphism spaces, ie exactly the functoriality map on $n$-cells:
\[
F^n: \End^n_{\LMod^n_R} (\LMod_R^{n-1}) \to \End^n_{\LMod^n_{\Dual^n R}} ({}_{\Dual^n R} Q),
\] 
This clearly gives a bimodule ${}_{\Dual^n R} Q_R$.

Applying this to our functor $F = {}_{\Dual^n R} k_R \otimes_R {-}$, we note that $F$ clearly corresponds to the bimodule given by Koszul duality, ${}_{\Dual^n R} k_R$. Also notice that
\[
\alpha' = F^n: \End^n_{\LMod^n_R} (\LMod_R^{n-1}) \to \End^n_{\LMod^n_{\Dual^n R}} (\LMod^{n-1}_{k, \aug}).
\]
This map can be identified with a map 
\begin{equation} \label{equation:unit}
R \to \Dual^n \Dual^n R     
\end{equation}
using our calculations in \ref{cor:morita1} and \ref{remark:morita2}.
This map MUST be the adjunct to the Koszul duality pairing 
\[
R \otimes \Dual^n R \to k,
\] since it must give $k$ the Koszul duality structure ${}_{\Dual^n R} k_R$ by what we said above. In other words, \ref{equation:unit} is the unit of the self-adjuntion of $\Dual^n$. Since $R$ is small, this unit is an equivalence (just combine \cite[Theorem~4.5.5]{DAGX} with \cite[Proposition~1.3.5]{DAGX}). Thus $\alpha'$, and thus $\alpha$, is an equivalence, as needed.

We've finally proved that $\theta$ is an equivalence. Finally we can conclude that $\Omega^n \ObjDefMap$ is an equivalence because $\Omega^n \ObjDefMap$ is induced by taking fibers of the vertical maps in the square
\[
\begin{tikzcd}
\End^n_{\C_R}(M_R) \ar[r, "\theta"] \ar[d] & \End^n_{\C_{\Dual^n R}}(M_{\aug}) \ar[d] \\
\End^n_{\C}(M) \ar[r, "\id"] & \End^n_{\C}(M).
\end{tikzcd}
\]
\end{proof}

\begin{remark}
Note that since we can identify $\theta$ and $\alpha$, our last square can be simplified to be
\[
\begin{tikzcd}
R \otimes \zeta(M) \ar[r, "\alpha"] \ar[d, "aug"] & \Dual^n \Dual^n R \otimes \zeta(M) \ar[d, "aug"] \\
\zeta(M) \ar[r, "\id"] & \zeta(M)
\end{tikzcd}
\]
which after taking fibers (which gives $\Omega^n \ObjDefMap)$, just says that $m_R \otimes \zeta(M) \simeq m_{\Dual^n \Dual^n R} \otimes \zeta(M)$.
\end{remark}

\begin{remark}
We can alternatively follow \cite[Proposition~5.3.19]{DAGX}, to prove this last step instead. We'll use the fiber sequence $m_A \to A \to k$ instead of the more restricted $k[n] \to k \oplus k[n] \to k$, to get the following:

We can identify the domain of $\Omega^n \ObjDefMap_{A}$ with
\[
\Maps_{\naugAlg{n}}(\Free^n(m_{A}^{\vee}), k \oplus \Z(M))
\]
and the codomain with 
\[
\Maps_{\naugAlg{n}}(\Dual^n(\Omega^n A), k \oplus \Z(M)).
\]
Using the functoriality of $\Omega^n \ObjDefMap_{A}$ in $(\C, M)$---the input for $\Z(M)$---and the Yoneda lemma, it follows that $\Omega^n \ObjDefMap_{A}$ is induced by a map 
\[
\Dual^n(\Omega^n A) \to \Free^n(k[m_{A}^{\vee}]).
\]
This is equivalent to having a map after passing to Koszul duals (for small algebras $A$) which gives 
\[
k \oplus m_{\Omega^n A} \to \Omega^n A.
\]
If one can show that this map is an inverse to the natural comparison map 
\[
\eta: \Omega^n A \to k \oplus m_{\Omega^n A}
\]
then we would be done. Here we are using the fact that $\Omega^n A$ is actually square zero, so $\eta$ is an equivalence.
\end{remark}
\subsection{Examples}
Here we list several example deformation problems that our theorem \ref{theorem:objdef2} characterizes. 

\begin{example}[Object in a 1-category]
Taking a $1$-category $\C \in \Prlk$ and an object $M \in \C$, we see that we recover Lurie's result about deforming an object in a category (\cite[Section~5.2]{DAGX}, \cite[Section~16.5]{SAG}). 

Thus our theorem \ref{theorem:objdef2} recovers many classical results, like deforming a quasicoherent module $M$ on a scheme $X$ over the dual numbers $k[\epsilon]$. We can use $\C = \QCoh(X)$, and we see that deformations of $M$ over $k[\epsilon]$ are given by maps 
\[
\Dual^{(1)}(k[\epsilon]) = k\langle \eta \rangle \to \End(M)
\]
where $k\langle \eta \rangle$ is the free associative algebra generated in cohomological degree $1$. Taking $\pi_0$ of the hom space 
\[
\Hom_k(k\langle \eta \rangle, \End(M))
\]
thus gives the first cohomology $H^1(\End(M))$, recovering the classical result. 
\end{example}

\begin{example}[Deformations of categories]
Taking the category to be $\Prlk$, we recover Lurie's result on deforming categories (\cite[Section~5.3]{DAGX}, \cite[Section~16.6]{SAG}). In particular, note that the center of a category $\D \in \Prlk$ is its Hochschild homology 
\[
HH(\D) = \End_{\End(\D)}(1_{\D}).
\]
\end{example}

\section{Simultaneous deformations}
Throughout this section, we assume we are given $\C \in \nPrlk{n}$ and an object $M \in \C$. We now shift our study to the situation of deforming an object $M$ and a category $\C$ together, considered as an object of $\nPrlkpoint{n}$. We aim to show such deformations are characterized by the $\E{n+1}$-algebra $\Z(\C, M)$, which can be described as the fiber of $\Z(\C) \to \Z(M)$ (see $\ref{eq:Z}$).

The argument again follows four steps, following ideas in \cite[Section~5.2, 5.3]{DAGX} and \cite[Section~4.1]{Blanc}: 
\begin{enumerate}
    \item Construct the functor $\SimDef_{(\C, M)}$.
    \item Prove $\SimDef_{(\C, M)}$ is $n + 1$-proximate.
    \item Construct the comparison map \[\SimDefMap: \SimDef_{(\C, M)} \to \Maps_{\naugAlg{n+1}}(\Dual^{n+1}({-}), k \oplus \Z(\C, M))\]
    \item Prove that $\SimDefMap$ is an equivalence.
\end{enumerate}
This idea is very similar to the case of deforming an object in an $n$-category and we hope to unify these two approaches in the future. 

\subsection{Constructing the functor \texorpdfstring{$\SimDef_{(\C, M)}$}{SimDef(C,M)}}

Like the $\ObjDef$ case, we start with the functor 
\[
\LMod^n: \Alg^{n+1}_k \longrightarrow \Alg(\nPrlk{n}) \simeq \Alg(\nPrlkpoint{n})
\]
given by $A \mapsto \LMod^n_A$, or $A \mapsto (\LMod^n_A, \LMod^{n-1}_A)$ in the pointed version (see \ref{definition:lmod}). 

\begin{definition}[Left Module categories]\label{definition:lcat}
Let $\LCatn$ be the pullback
\[
\begin{tikzcd}
\LCatn \ar[r] \ar[d] \arrow[dr, phantom, "\scalebox{1.5}{$\lrcorner$}" , very near start, color=black]
& \LMod(\nPrlk{n}) \ar[d]
\\ \nAlg{n+1} \ar[r] & \Alg(\nPrlk{n})
\end{tikzcd}
\]
and let $\LCatnstar$ be the pullback
\[
\begin{tikzcd}
\LCatnstar \ar[r] \ar[d] \arrow[dr, phantom, "\scalebox{1.5}{$\lrcorner$}" , very near start, color=black]
& \LMod(\nPrlkpoint{n}) \ar[d]
\\ \nAlg{n+1} \ar[r] & \Alg(\nPrlkpoint{n})
\end{tikzcd},
\]
where the left leg is given by the usual cocartesian fibration \cite[Definition~4.2.1.13]{HA}. There's a obvious projection $\LCatnstar \to \LCatn$ that forgets the basepoint.
\end{definition}
\begin{remark}
Intuitively, objects of $\LCatnstar$ consists of $4$-tuples $(A, \D, N, \eta)$ where $N \in \D$, $A \in \nAlg{n+1}$, and $\eta$ is a left action of $\LMod^{n}_A$ on $(\D, N) \in \nPrlkpoint{n}$. 
\end{remark}

\begin{remark}
Dually, we can use right modules instead by replacing the right vertical leg by
\[
\RMod(\nPrlkpoint{n}) \longrightarrow \Alg(\nPrlkpoint{n})
\] 
So we can analoguously define $\RCatnstar$. This has objects $(A, \D, M, \eta)$ where $N \in \D$, $A \in \nAlg{n+1}$, and $\eta$ is a {\it right} action of $\LMod^{n}_A$ on $(\D, N) \in \nPrlkpoint{n}$. 

Similarly we can analogously define $\RCatn$.
\end{remark}

First, just like with ObjDef (\ref{definition:objfib}), we first define the associated left fibration to SimDef. 

\begin{definition}[Simultaneous deformation fibration]\label{defintion:sDef}
Let $\LCatnstar \to \nAlg{n+1}$ be the cocartesian fibration defined above (\ref{definition:lcat}). We can get a left fibration by restricting to only cocartesian arrows:
\begin{equation}\label{equation:lcatfib}
\LCatnstarcocart \longrightarrow \nAlg{n+1}.
\end{equation}
Finally, let $\sDef[\C, M]$ be the slice
\[
(\LCatnstarcocart)_{/(k, \C, M)}.
\]
By slicing the fibration (\ref{equation:lcatfib}), we see that we have a left fibration
\begin{equation}\label{equation:sDeffib}
\sDef[\C, M] \to \nAlgslice{n+1} \simeq \naugAlg{n+1}.
\end{equation}
\end{definition}

Now we are ready to construct $\SimDef_{(\C, M)}$.
\begin{construction}[$\SimDef$] 
We look at our given pair $(\C, M)$, with natural $\LMod^{n}_k$ action, which we'll denote as $(k, \C, M) \in \LCatnstar$, the action being implicit. 

By straightening, the left fibration (\ref{equation:sDeffib})
\[
\sDef[\C, M] \to \naugAlg{n+1}
\]
classifies a functor $\widetilde{\sDef}[\C, M]: \naugAlg{n+1} \to \bigSpc$. Here $\bigSpc$ is the category of not-necessarily $U_0$-small spaces. Finally, by restricting to small algebras, we get the functor we wanted: 
\[
\SimDef_{(\C, M)}: \nsmAlg{n+1} \to \bigSpc.
\]
\end{construction}
\begin{remark}
Notice that given an algebra $A \in \naugAlg{n+1}$, 
\[
\SimDef_{(\C, M)}(A) \simeq \nPrlpoint{n}{A} \times_{\nPrlkpoint{n}} \{(\C, M)\},
\] 
as mentioned in the introduction.
\end{remark}
\subsection{Proving \texorpdfstring{$\SimDef_{(\C, M)}$}{SimDef(C,M)} is \texorpdfstring{$n+1$}{n+1}-proximate}

To prove this statement, we show the existence of a fiber sequence of deformation functors, generalizing Proposition 4.3 of \cite{Blanc}.

We begin with constructing the maps. 

\begin{construction}[Comparison with ObjDef]
We construct the projection 
\[
\SimDef_{(\C, M)} \longrightarrow \ObjDef_{\C \in \nPrlk{n}}.
\]
which intuitively just forgets the "point" $M$ and its deformation, ie it sends a simultaneous deformation $(A, \C_A, M_A)$ and forgets the $M$-deformation $M_A$. 

More precisely, first we use the projection 
\[
\LCatnstar \longrightarrow \LCatn
\]
which commutes to the projections to $\nAlg{n+1}$, which induces a projection of slices
\[
\sDef[\C, M] \longrightarrow \Deform[\C \in \nPrlk{n}]
\]
over $\naugAlg{n+1}$. This induces a projection
\begin{equation}\label{map:fiber1}
\SimDef_{(\C, M)} \longrightarrow \ObjDef_{\C \in \nPrlk{n}}.
\end{equation}
of functors $\nsmAlg{n+1}\to \bigSpc$
\end{construction}

\begin{remark}
The codomain of this projection, in the case that $n=1$, is usually called $\CatDef_{\C}$. 
\end{remark}

Next we analyze the fiber of this projection.
\begin{construction}
The fiber of the projection
\[
\nPrlkpoint{n} \longrightarrow \nPrlk{n}
\]
at the category $\C$ is equivalent to the Kan complex $\Maps_{\nPrlk{n}}(\Mod^n_k, \C) \simeq \C^{\simeq}$. The fiber map $\C^{\simeq} \to \nPrlkpoint{n}$ sends $M \in \C$ to $(\C, M)$. This induces a functor 
\[
\Deform[M \in \C] \longrightarrow \sDef[\C, M]
\]
commuting with the projections to $\naugAlg{n+1}$.
This induces a natural transformation 
\begin{equation}\label{map:fiber2}
\ObjDef_{M \in \C} \longrightarrow \SimDef_{(\C, M)} 
\end{equation}
\end{construction}
By our construction, it is directly obvious that we have a fiber sequence:

\begin{proposition}\label{prop:fiber}
The natural transformations constructed above (\ref{map:fiber1}, \ref{map:fiber2}) fit into a fiber sequence of functors $\nsmAlg{n+1} \to \bigSpc$:
\[
\ObjDef_{M \in \C} \to \SimDef_{(\C, M)} \to \ObjDef_{\C \in \nPrlk{n}}.
\]
\end{proposition}

Now we show some consequences of having this fiber sequences. First, it proves what we wanted to show:
\begin{proposition}
$\SimDef_{(\C, M)}$ is an $n+1$-proximate fmp.
\end{proposition}
\begin{proof}
Since $\ObjDef_{M \in \C}$ is $n$-proximate and $\ObjDef_{\C \in \nPrlk{n}}$ is $n+1$-proximate (\ref{section:2.2}), we see that since taking loop spaces and pullbacks preserve limits, that $\SimDef_{(\C, M)}$ must also be $n+1$-proximate. 
\end{proof}

Since the completion functor from $n+1$-proximate fmps to fmps is limit-preserving, we can also easily see

\begin{proposition}
The fiber sequence we constructed descends to fmp completions, and we have a natural comparison of fiber sequences:
\[
\begin{tikzcd}
\ObjDef_{M \in \C} \ar[r] \ar[d] & \SimDef_{(\C, M)} \ar[r] \ar[d] & \ObjDef_{\C \in \nPrlk{n}} \ar[d]
\\ \fmpObjDef^{\E{n+1}}_M \ar[r] & \fmpSimDef_{(\C, M)} \ar[r] & \fmpObjDef_{\C \in \nPrlk{n}}
\end{tikzcd}
\]
where the vertical maps are the units for the fmp-completion functor. 
\end{proposition}
\subsection{Constructing the comparison map \texorpdfstring{$\SimDefMap$}{b}}
For this section, we'll use the following notation:

\begin{notation}
Given an $n$-category $\C$ with a left action by an $\E{n}$-algebra $A$ (in other words, a left action by $\LMod^n_A)$ and a right action by an $\E{n}$-algebra $B$, we emphasize this structure as so:
\[
{}_A [\C]_B
\]
\end{notation}

Let $\lambda^{n+1}: \pair^{n+1} \to \naugAlg{n+1} \times \naugAlg{n+1}$ be the pairing of categories inducing $\E{n+1}$-Koszul duality (\cite[Construction~5.2.5.32]{HA}). Objects of $\pair^n$ intuitively consist of two algebras $A, B \in \naugAlg{n+1}$ along with an augmentation of their tensor product: $A \otimes_k B \to k$. This gives $k$ the structure of a $A \otimes_k B$ module.

Let $A \otimes_k B \to k$ be an $\E{n+1}$-pairing between $A, B \in \nsmAlg{n+1}$ (so it is an object of $\pair^{n+1}$). Suppose we are given a simultaneous deformation $(A, \C_A, M_A) \in \SimDef_{(\C, M)}(A)$ (where the actions and augmentation equivalences are suppressed).
Notice that 
\[
(\C_A, M_A) \otimes \LMod^{n}_B \in {}_{\LMod^{n}_A \otimes \LMod^{n}_B} \BiMod_{\LMod^{n}_B}(\nPrlkpoint{n})
\]
Here we let $\LMod^{n}_B$ stand in for the pointed category $(\LMod^{n}_B, \LMod^{n-1}_B) \in \nPrlkpoint{n}$.

Thus we have 
\[
(\C, M) \simeq \LMod^{n}_k \underset{A \otimes B}{\otimes} [(\C_A, M_A) \otimes \LMod^{n}_B] \in \RMod_{\LMod^{n}_B}(\nPrlkpoint{n})
\]
where the equivalence uses the given pairing and the augmentation equivalences. This construction gives a right $\LMod^n_B$ action on $(\C, M)$. Notice we gave $\Mod^n_k$ a right $\LMod^{n}_A \otimes \LMod^{n}_B$ module structure using the augmentation. 

In fact we have just a little more, we know that the action of $\LMod^n_B$ on $M$ is "trivial": $(\C_A, M_A)$ can be written as a map
\[
{}_{A} [\LMod^n_A]_k \longrightarrow {}_{A} [\C_A]_k
\]
The shorthand subscripts on the left denote left actions of $\LMod^n_A$ and analogously, the right subscripts show right actions.

Then the above process can be seen as tensoring on the left by $\Mod^n_k$, this time seen as having a left $\LMod^n_{B^{\op}}$ and right $\LMod^n_A$ action (which is equivalent to the right action of $A \otimes B$. All this distinction of $\op$-algebras only matters for the very trivial case $n = 0$). Thus we get
\[
 {}_{B^{\op}}[{\Mod^n_k}]_{A} \underset{A}{\otimes}  {}_{A}[\LMod^n_A]_k \longrightarrow  {}_{B^{\op}}[{\Mod^n_k}]_{A} \underset{A}{\otimes} {}_{A}[\C_A]_k
\]
where our tensors are over $A$. This reduces to 
\[
{}_{B^{\op}}[{\Mod^n_k}]_k \longrightarrow {}_{B^{\op}} [\C]_k
\]
which shows the right $B$-action (or left $B^{\op}$-action). Notice that the action of $B$ on $\Mod^n_k$ is "trivial" since it is no longer coupled with the action of $A$. Hence it is image, which points out $M$, also has trivial action in this way. 

This construction produces a functor:
\begin{equation}\label{eq:pairing}
\sDef[\C, M] \times \pair^{n+1} \to \sDef[\C, M] \times (\RCatnstartriv \times_{\nPrlk{n}} \{(\C,M)\}) 
\end{equation}
where $\RCatnstartriv$ is the category of pointed $n$-categories with right actions by {\it augmented algebras} which are trivial on the given point. Each object is an object of $\RCatnstar$ with extra triviality data. In other words, $(\D, E)$ with $R$-action (where $R$ is augmented) is trivial when the map 
\[
\Mod^n_k \longrightarrow \D
\]
picking out $E$ factors through the right $R$-module map
\[
[\Mod^n_k]_R \longrightarrow [\D]_R,
\]
where the action of $R$ on $\Mod^n_k$ is given by the augmentation map. The explicit definition is given as follows:
\begin{construction}[$\RCatnstartriv$]
We have a functor 
\[
\Triv: \naugAlg{n+1} \to \RMod(\nPrlk{n})
\]
sending an augmented algebra $A$ to the category with right augmentation $A$-action $\LMod_k$. Let 
\[
\pi: \RMod(\nPrlk{n}) \to \Alg(\nPrlk{n})
\] be the natural projection. Finally, we have two maps from $\naugAlg{n+1}$. First we have
\[
\aug: \naugAlg{n+1} \to \Arr(\RMod(\nPrlk{n}))
\]
taking $A$ to the augmentation functor $\LMod_A \to \LMod_k$, equipped with natural right $A$-action and right augmentation $A$-action respectively. Secondly we have
\[
d: \naugAlg{n+1} \to \Fun(\Delta^2, \Alg(\nPrlk{n}))
\]
sending $A$ to the degenerate triangle
\[
\begin{tikzcd}
& \LMod_A \ar[rd, "="] \\
\LMod_A \ar[ru, "="] \ar[rr, "="] & & \LMod_A
\end{tikzcd}
\]
Finally we can define:
Let $\RCatnstartriv$ be the pullback:
\[
\begin{tikzcd}
\RCatnstartriv \ar[r] \ar[d] \arrow[dr, phantom, "\scalebox{1.5}{$\lrcorner$}" , very near start, color=black] & \Fun(\Delta^2, \RMod(\nPrlk{n})) \ar[d, "{(\ev_{[0, 1]}, \pi)}"]
\\ \naugAlg{n+1} \ar[r, "{(\aug, d)}"] & {\Arr(\RMod(\nPrlk{n})) \times \Fun(\Delta^2, \Alg(\nPrlk{n}))}
\end{tikzcd}
\]
Of course, we can analogously define the version with left actions instead.
\end{construction}

Now continuing our construction, the functor we constructed (\ref{eq:pairing}) is a left representable pairing of categories, which induces a duality functor:
\begin{equation} \label{map:dual}
\Dual^{n+1}: \sDef[\C, M]^{\op} \to \RCatnstartriv \times_{\nPrlkpoint{n}} \{(\C, M)\} 
\end{equation}

Notice the codomain of this functor is equivalent to $\naugAlgslice{n+1}{\Z(\C, M)}$, where $\Z(\C, M)$ is as defined in \ref{eq:Z}.  

We have a square:

\[
\begin{tikzcd}
\sDef[\C, M]^{\op} \ar[r, "\Dual^{n+1}"] \ar[d] & \naugAlgslice{n+1}{\Z(\C, M)} \ar[d]
\\ \naugAlg{n+1}^{\op} \ar[r, "\Dual^{n+1}"] & \naugAlg{n+1}
\end{tikzcd}
\]
Here the top horizontal functor is the duality functor we just defined (\ref{map:dual}), the bottom functor is $\E{n+1}$-Koszul duality functor. The left and right vertical maps are canonical Cartesian fibrations. By restricting to small algebras and using straightening/unstraightening, we finally get a comparison morphism of the two induced functors $\nsmAlg{n+1} \to \bigSpc$:
\begin{equation} \label{map:simbeta}
\SimDefMap: \SimDef_{(\C, M)} \to \Maps_{\naugAlg{n+1}}(\Dual^{n+1}({-}), k \oplus \Z(\C, M)) 
\end{equation}

\subsection{Proving \texorpdfstring{$\SimDefMap$}{bSim} induces an equivalence}
Now we can finally prove our main theorem:
\begin{theorem} \label{theorem:simdef2}
The map $\SimDefMap$ (\ref{map:simbeta}) induces an equivalence 
\[
\fmpSimDef_{(\C, M)} \to \Maps_{\naugAlg{n+1}}(\Dual^{n+1}({-}), k \oplus \Z(\C, M))
\]
of formal moduli problems.
\end{theorem}
We'll use some notation here:
\begin{notation}[(Fmp associated to an algebra)]
Given an $\E{n}$ augmented algebra $R$, let $\Psi^n_R$ denote the $\E{n}$-fmp associated to $R$. In other words, 
\[
\Psi^n_R := \Maps_{\naugAlg{n}}(\Dual^{n}{-},R).
\]
\end{notation}

To do this, we only have to show that after taking $n+1$-fold loop spaces, we have an equivalence
\[
\Omega^{n+1} \SimDef_{(\C, M)} \to \Omega^{n+1} \Psi^{n+1}_{k \oplus \Z(\C, M))}
\] 
since $\SimDef_{(\C, M)}$ is an $n$-proximate fmp. We use Lurie's Proposition 1.2.10 in \cite{DAGX} to reduce to the cases where the algebra is $k \oplus k[m]$, as values on these algebras determine the tangent complex in our current context. 

So we've reduced the theorem to the following proposition:
\begin{proposition}
Let $A := k[m]$ and $B := k \oplus k[m+n+1]$. Then $A$ is the $n+1$-th loop space of $B$. This induces the following diagram: 
\[
\begin{tikzcd}\label{diagram:simdef}
\SimDef_{(\C, M)}(A) \ar[r, "\SimDefMap"] \ar[d] \arrow[dr, phantom, "\scalebox{1.5}{$\lrcorner$}" , very near start, color=black] & \Psi^{n+1}_{k \oplus \Z(\C, M))}(A) \ar[d, "\sim"]
\\ \Omega^{n+1} \SimDef_{(\C, M)}(B) \ar[r, "\Omega^{n+1} \SimDefMap"] & \Omega^{n+1} \Psi^{n+1}_{k \oplus \Z(\C, M))}(B)
\end{tikzcd}
\]
Then the bottom map of this diagram is an equivalence.
\end{proposition}

\begin{proof}
We have the fiber sequence $k[m+n+1] \to B \to k$. Let $\C_B := \LMod^n_B(\C)$ and $M_B := \LMod^{n-1}_B \otimes M$. 

Notice that $\Omega^{n+1} \SimDef_{(\C, M)}(B)$ can be identified with the fiber of the augmentation map
\[
\End^{n+1}_{\nPrlpoint{n}{B}}(\C_B, M_B) \to \End^{n+1}_{\nPrlkpoint{n}}(\C, M).
\]
We'll suppress the $\nPrlkpoint{n}$ and $\nPrlpoint{n}{B}$ from now on. 
These two objects are both fibers:
\[
\End^{n+1}(\C_B, M_B) \simeq 
\begin{tikzcd}
\fib(\End^{n+1}(C_B) \ar[r, "\ev_{M_B}"] &\End^n_{\C_B}(M_B))
\end{tikzcd}
\]
\[
\End^{n+1}(\C, M) \simeq 
\begin{tikzcd}
\fib(\End^{n+1}(C) \ar[r, "\ev_M"] &\End^n_{\C}(M))
\end{tikzcd}
\]
where we suppressed the $\nPrlk{n}$ in the subscripts.
Therefore, $\Omega^{n+1} \SimDef_{(\C, M)}(B)$ can be identified with the fiber of the comparison map of the square
\begin{equation}\label{diagram:end}
\begin{tikzcd}
\End^{n+1}(C_B) \ar[r, "\ev_{M_B}"] \ar[d, "aug"] &\End^{n}_{C_B}(M_B) \ar[d, "aug"] \\ 
\End^n(\C_B) \ar[r, "\ev_M"] &\End^{n+1}_{\C}(M)
\end{tikzcd}
\end{equation}
where by "comparison map" we mean the map from the top left corner to the pullback of the bottom and right legs.

Notice that if we take the fiber of the left map in this square, we get 
\[\Omega^{n+1} \ObjDef_{\C \in \nPrlk{n}}(B)\]
and if we take the fiber of the right map in this square, we get
\[
\Omega^n \ObjDef_M(B).
\]
This shows we have a fiber sequence
\begin{equation}\label{equation:fiber}
\Omega^{n+1} \SimDef_{(\C, M)}(B) \to \Omega^{n+1} \ObjDef_{\C \in \nPrlk{n}}(B) \to \Omega^n \ObjDef_M(B).  
\end{equation}
which is natural in $B$, or more generally, any $\E{n+1}$ algebra. Notice this sequence is closely related to the sequence we used in \ref{prop:fiber}. The second map is induced by evaluation at $M$, just as in the square \ref{diagram:end}. 

Let $\Psi^n_A$ denote the $\E{n}$-fmp associated to the augmented algebra $A$. In other words, 
\[
\Psi^n_A := \Maps_{\naugAlg{n+1}}(\Dual^{n+1}(B),A)
\]

Clearly by using the various comparison maps for each fmp (\ref{map:objbeta}, \ref{map:simbeta}), we get the following diagram
\begin{equation}\label{diagram:rectangle}
\begin{tikzcd}
\Omega^{n+1} \SimDef_{(\C, M)}(B) \ar[d, "\Omega^{n+1} \SimDefMap"] \ar[r] & \Omega^{n+1} \ObjDef_{\C \in \nPrlk{n}}(B) \ar[d, "\sim"] \ar[r] & \Omega^n \ObjDef_M(B) \ar[d, "\sim"]\\
\Omega^{n+1} \Psi^{n+1}_{k \oplus \Z(\C, M)}(B) \ar[r] & \Omega^{n+1} \Psi^{n+1}_{k \oplus \Z(\C)}(B) \ar[r]
& \Omega^{n} \Psi^{n}_{k \oplus \Z(M)}(B)
\end{tikzcd}
\end{equation}

If we can show that the bottom sequence is also a fiber sequence, we'll be done, as the comparison $\Omega^{n+1} \SimDefMap$ on the left would have to be an equivalence. 

Since the $\Psi$'s are fmps, we can push in the loop spaces to act on $B$, giving us this equivalent bottom sequence:
\[
\Psi^{n+1}_{k \oplus \Z(\C, M)}(A) \to \Psi^{n+1}_{k \oplus \Z(\C)}(A) \to \Psi^{n}_{k \oplus \Z(M)}(k \oplus k[m+1])
\]
Next we can combine the identify $\Dual^{n+1}(k \oplus k[N])$ with $\Free^{n+1}(k[-N -n -1])$ with the free-forgetful adjunction (with nonunital algebras) to identify the sequence with 
\[
\Maps_k(k[-m-n-1], \Z(\C, M)) \to \Maps_k(k[-m-n-1], \Z(\C)) \to \Maps_k(k[-m-n-1], \Z(M))
\]
where the second map is by evaluation at $M$, just like in \ref{diagram:end}. This is clearly a fiber sequence since 
\[
\Z(\C, M) \to \Z(\C) \to \Z(M)
\]
is a fiber sequence in nonunital algebras (where the second map is evaluation at $M$).
\end{proof}

\subsection{Examples} \label{section:3.5}
Now we give an application of theorem \ref{theorem:simdef2}.
\begin{example}[Monoidal categories and algebras]
Our main example is deforming $\E{n}$-monoidal $m$-categories, which includes the case of deforming $\E{n}$-algebras as the $m=0$ case. These are very important in the study of shifted sympletic structures. For example, \cite{PTVV} discusses the relation between $n$-Poisson structures on derived affine stacks and the deformation theory of $\E{n}$-monoidal categories. They are also important in the study of Quantization, as discussed in \cite{Toen}. 

Given an $\E{n}$-monoidal $m$-category $\D^{\otimes}$, its deformation theory is the same as the simulatenous deformation theory of the pair $(\LMod^n_{\D}, \LMod^{n-1}_{\D})$. Namely, theorem \ref{theorem:simdef2} says that $\E{n}$-monoidal deformations of the $\E{n}$-monoidal $m$-category $\D^{\otimes}$ are characterized by 
\[
\Z(\LMod^n_{\D}, \LMod^{n-1}_{\D}) = \fib(\Z(\LMod^n_{\D}) \to \Z(\LMod^{n-1}_{\D})).
\]
In the $\E{1}$-algebra case, we see that deformations of an algebra $A$ are characterized by 
\[
\Z(\LMod_A, A) = \fib(HH(A) \to A)
\]
where $HH(A)$ denotes the Hochschild complex of $A$ (ie the derived center of $\LMod_A$). We can recover the classical result that first-order deformations of $A$ are characterized by $HH^2(A)$ when $A$ is connective. See for example \cite{Fox} for the case that $A$ is concentrated in degree $0$. Note that in this case, 
\[
HH^2(A) = \pi_0(\Hom_k(k[-2], HH(A))) \simeq \pi_0(\Hom_k(k[-2], \fib(HH(A) \to A)))
\]
since $A$ is connective, by using the long exact sequence of cohomology groups. Then we use our theorem \ref{theorem:simdef2} to see:
\begin{align*}
\fmpAlgDef_A(k \oplus k[0]) & \simeq \Hom_{\naugAlg{2}}(\Free(k[-2]), \fib(HH(A) \to A))\\
& \simeq \Hom_k(k[-2], \fib(HH(A) \to A)),
\end{align*}
and hence
\[
\pi_0(\fmpAlgDef_A(k \oplus k[0])) \simeq HH^2(A).
\]
We sketch the argument for deforming $\E{n}$-monoidal $m$-categories. Given an $\E{n}$-monoidal $m$-category $\D^{\otimes}$ (note that we drop the superscript sometimes, especially when we need $\D^\otimes$ itself as a subscript), we let $\nMonDef{n}_\D$ be the functor that assigns to each small $\E{n+m+1}$-algebra $B$ to the groupoid core of
\[
\nMonCat{n}_B \times_{\nMonCat{n}_k} \{\D^\otimes\},
\] 
where $\nMonCat{n}_B$ denotes the category of $\E{n}$-monoidal $m$-categories with a central $B$-action, and the map $\nMonCat{n}_B \to \nMonCat{n}_k$ is given by the augmentation map $B \to k$. 
Observe that $\LMod^n_{\D}$ is now an $n+m$-category as it is a category consisting of $n+m-1$-categories with a $\D$-central action. Hence notice that taking $n$-fold endomorphism spaces in this category produces $m$-categories.
Notice one can define using a left fibration like we do with $\ObjDef$ (\ref{construction:objdef}). 

The key insight is that we have a comparison 
\[
\nMonDef{n}_{\D} \to \ObjDef_{\LMod^n_{\D} \in \nPrlk{n+m}} =: \CatDef_{\LMod^n_{\D}}
\]
given by taking a deformation $\D^\otimes_B \in \nMonDef{n}_{\D}(B)$ to $\LMod^n_{\D_B}$. Further there's a map 
\[
\ObjDef_{\LMod^{n-1}_{\D} \in \LMod^n_{\D}} \to \nMonDef{n}_{\D}
\]
which sends a deformation $M \in B \otimes \LMod^n_\D$ to $\End^n_{\D}(M)$. This is an $m$-category by what we said above. Notice that the composition
\[
\ObjDef_{\LMod^{n-1}_{\D} \in \LMod^n_{\D}} \to \nMonDef{n}_{\D} \to \CatDef_{\LMod^n_{\D}}
\]
is a fibration because the monoidal categories $\D^\otimes_B$ that get sent to the basepoint $\LMod^n_{B \otimes \D}$ of $\CatDef_{\LMod^n_{\D}}$ are exactly characterized by objects $M \in \ObjDef_{\LMod^{n-1}_{\D}}(B)$ by Morita-equivalence arguments.

One can construct a comparison map 
\[
\nMonDef{n}_{\D} \to \SimDef_{(\LMod^n_{\D}, \LMod^{n-1}_{\D})}
\]
in the following way: map $\D_B \in \nMonDef{n}_{\D}(B)$ to the pair 
\[
(\LMod^n_{\D_B}, \LMod^{n-1}_{\D_B}) \in \SimDef_{(\LMod^n_{\D}, \LMod^{n-1}_{\D})}(B).
\]
Then we have a comparison of fiber sequences
\[
\begin{tikzcd}
\ObjDef_{\LMod^{n-1}_{\D} \in \LMod^n_{\D}} \ar[r] \ar[d] & \nMonDef{n}_{\D} \ar[r] \ar[d] & \CatDef_{\LMod^n_{\D}} \ar[d]
\\ \ObjDef_{\LMod^{n-1}_{\D} \in \LMod^n_{\D}} \ar[r] & \SimDef_{(\LMod^n_{\D}, \LMod^{n-1}_{\D})} \ar[r] & \CatDef_{\LMod^n_{\D}}.
\end{tikzcd}
\]
This is an equivalence after passing to fmp completions, as the left and right legs are equivalences.
\end{example}

\DeclareFieldFormat{labelnumberwidth}{{#1\adddot\midsentence}}
\printbibliography
\end{document}